\documentclass{article}

\usepackage{latexsym,theorem,graphics,xspace,ifthen}
\usepackage{amssymb,amsmath,amscd}

\def\today{\ifcase\month\or
 January\or February\or March\or April\or May\or June\or
 July\or August\or September\or October\or November\or December\fi
 \space\number\day, \number\year}

\numberwithin{equation}{section}
\DeclareFontFamily{OT1}{rsfs}{}
\DeclareFontShape{OT1}{rsfs}{m}{n}{ <-7> rsfs5 <7-10> rsfs7 <10-> rsfs10}{} 
\DeclareMathAlphabet{\mathscr}{OT1}{rsfs}{m}{n}

\newcommand{\brckt}[1]{\ifthenelse{\equal{#1}{}}{}{(#1)}}

\newcommand{\db}[1][{}]{\ensuremath{\bar{\partial}^{#1}_{b}}\xspace}
\newcommand{\dbs}[1][{}]{%
\ifthenelse{\equal{#1}{}}{\ensuremath{\bar{\partial}^*_{b}}\xspace}
{\ensuremath{(\bar{\partial}^{#1}_b)^*}\xspace}}%

\newcommand{\vtb}{\ensuremath{\vartheta_b}\xspace}

\newcommand{\dist}[3][U]{\ensuremath{\text{dist}_{#1}\left(#2,#3\right)}\xspace}

\newcommand{\dD}[1][{}]{{\partial_{#1} D}}

\newcommand{\ddDn}[1][\cdot]{\dist[]{\cdot}{\partial D}}

\newcommand{\norm}[3][{}]{\ensuremath{\big\|#2\big\|_{#3}^{#1}}\xspace}
\newcommand{\snorm}[3][{}]{\ensuremath{\big|#2\big|_{#3}^{#1}}\xspace}
\newcommand{\normV}[3][{}]{\ensuremath{\left\|#2\right\|_{#3}^{#1}}\xspace}
\newcommand{\snormV}[3][{}]{\ensuremath{\left|#2\right|_{#3}^{#1}}\xspace}

\newcommand{\crT}[1][{}]{\ensuremath{{^o T_{#1}^{\prime}}}\xspace}
\newcommand{\acrT}[1][{}]{\ensuremath{{^o T_{#1}^{\prime \prime}}}\xspace}

\newcommand{\dom}[1]{\ensuremath{\text{Dom}(#1)}\xspace}
\newcommand{\domp}[1]{\ensuremath{\text{Dom}_+(#1)}\xspace}

\newcommand{\Ker}[1]{\text{Ker}(#1)}
\newcommand{\rng}[1]{\text{Range}(#1)}

\newcommand{\real}[1]{\text{Re}\left(#1 \right)}
\newcommand{\imag}[1]{\text{Im}\left(#1 \right)}

\newcommand{\nW}{\ensuremath{\overline{W}}\xspace}

\newcommand{\bt}[1]{{\bar{#1}}}
\newcommand{\Bt}[1]{{\overline{#1}}}

\newcommand{\za}[1][\s]{{\ua(#1)}}

\newcommand{\ua}{{\alpha}}
\newcommand{\ub}{{\beta}}

\newcommand{\ud}{{\delta}}
\newcommand{\ul}{{\lambda}}

\newcommand{\ut}{{\theta}}

\newcommand{\vs}{\ensuremath{{\varsigma}}\xspace}

\newcommand{\e}{{\epsilon}}

\newcommand{\dO}{\ensuremath{{\partial \Omega}}\xspace}

\newcommand{\Pt}[1][{\s}]{\ensuremath{{P}^\top_{#1}}\xspace}
\newcommand{\Pb}[1][{\s}]{\ensuremath{{P}^\bot_{#1}}\xspace}

\newcommand{\one}{\ensuremath{\mathsf{1}}\xspace}

\newcommand{\bQ}[2][{}]{%
\ifthenelse{\equal{#1}{}}{\ensuremath{\mathcal{Q}_{#2}}\xspace}
{\ensuremath{\mathcal{Q}_{#2}\big( #1 \big)}\xspace}}%

\newcommand{\bI}[1][{}]{%
\ifthenelse{\equal{#1}{}}{\ensuremath{\mathcal{I}}\xspace}{\ensuremath{\mathcal{I}\big(#1\big)}\xspace}}%

\newcommand{\btQ}[2][{}]{%
\ifthenelse{\equal{#1}{}}{\ensuremath{\Bt{\mathcal{Q}}_{#2}}\xspace}
{\ensuremath{\Bt{\mathcal{Q}}_{#2}\big( #1 \big)}\xspace}}%

\newcommand{\tQ}[2][{}]{\ensuremath{\widetilde{Q}^{#1}_{K} \big( #2 \big)}\xspace}

\newcommand{\CiD}{\ensuremath{C^\infty(\overline{D})}\xspace}
\newcommand{\CpD}{\ensuremath{C^\infty_+(\overline{D})}\xspace}
\newcommand{\CvD}{\ensuremath{C^\infty_v(\overline{D})}\xspace}

\newcommand{\Ci}[1]{\ensuremath{C^\infty(\overline{#1})}\xspace}
\newcommand{\Cv}[1]{\ensuremath{C^\infty_v(\overline{#1})}\xspace}
\newcommand{\Cic}[1]{\ensuremath{C^\infty_0\brckt{#1}}\xspace}

\newcommand{\CpO}[1][q]{\ensuremath{C^{(0,#1)}_+(\overline{\Omega})}\xspace}

\newcommand{\ip}[3]{\ensuremath{\big( \, {#1} \, , \, {#2} \, \big)_{#3}}\xspace}

\newcommand{\aip}[3]{\ensuremath{\langle \, {#1} \, , \, {#2} \, \rangle_{#3}}\xspace}

\newcommand{\boxb}[1][{}]{\ensuremath{\square^{#1}_b}\xspace}
\newcommand{\boxbt}{\ensuremath{\Bt{\square}_b}\xspace}

\newcommand{\upp}{\ensuremath{^{\prime}}\xspace}
\newcommand{\dupp}{\ensuremath{^{\prime \prime}}\xspace}

\newcommand{\rn}[1]{\ensuremath{\mathbb{R}^{#1}}\xspace}

\newcommand{\sn}[1]{\ensuremath{{\mathbb{S}^{#1}}}\xspace}
\newcommand{\cn}[1]{\ensuremath{\mathbb{C}^{ #1}}\xspace}
\newcommand{\hn}[1]{\ensuremath{\mathbb{H}^{ #1}}\xspace}
\newcommand{\hy}[1]{\ensuremath{\mathbb{U}^{ }}\xspace}

\newcommand{\LMain}[2]{\ensuremath{L^2_{#1}#2}\xspace}

\newcommand{\lO}[1][\ut]{\LMain{#1}{(\Omega)}}

\newcommand{\lOa}[1]{\LMain{\tLa}{}}

\newcommand{\xX}{\ensuremath{\mathcal{X}}\xspace}

\newcommand{\NabH}[1][H]{\ensuremath{\nabla_{ [#1]}}\xspace}

\renewcommand{\theenumi}{\alph{enumi}}

\newcommand{\WjMain}[3]{\ensuremath{\mathscr{W}^{#1}_{#2}#3}\xspace}

\newcommand{\Wj}[2][D]{\WjMain{#2}{}{\brckt{#1}}}
\newcommand{\Wpj}[2][D]{\WjMain{#2}{+}{\brckt{#1}}}

\newcommand{\WjR}[1]{\WjMain{#1}{}{}}
\newcommand{\WjRD}[1]{\WjMain{#1}{}{(D)}}
\newcommand{\WjsD}[1]{\WjMain{#1}{\s}{(D)}}
\newcommand{\Wjs}[2][{}]{\WjMain{#2}{\s}{#1}}
\newcommand{\WlD}[1]{\WjMain{#1,loc}{\s}{(\Bt{D})}}

\newcommand{\WoX}[3]{\ensuremath{\mathring{\mathscr{W}}_{#3}^{#1}#2}\xspace}
\newcommand{\Wo}[1][1]{\ifthenelse{\equal{#1}{1}}%
{\WoX{#1}{(D)}{}}%
{\WoX{#1}{(D)}{\s}}}%
\newcommand{\Won}[2]{\WoX{#1}{\brckt{#2}}{}}
\newcommand{\ED}[1]{\ensuremath{E\brckt{#1}}\xspace}
\newcommand{\V}[1][q]{\ensuremath{\mathscr{V}^{#1}}\xspace}
\newcommand{\Vv}[1][q]{\ensuremath{\mathscr{V}_{\nu}^{#1}}\xspace}

\newcommand{\s}{{\sigma}}
\newcommand{\mE}[1][\s]{\mathcal{E}\brckt{#1}}
\newcommand{\Lk}[1]{\ensuremath{L_{(#1)}}\xspace}
\newcommand{\Gm}[1][\s]{\Gamma\brckt{#1}}
\newcommand{\Go}[1][q]{\gamma_0^{#1}}
\newcommand{\G}[1][\s]{G(#1)}
\newcommand{\uls}[1][\s]{\ul(#1)}
\newcommand{\Ws}{W_\s}
\newcommand{\Wa}{W_{\ua}}
\newcommand{\nWs}[1][\s]{\Bt{W}_{\!#1}}
\newcommand{\ulsq}[1][\s]{|\uls[#1]|^{1/2}}
\newcommand{\R}[2][\s]{\ensuremath{\mathscr{R}^{#2}_{#1}}\xspace}
\newcommand{\Rp}[1]{\ensuremath{\mathscr{R}^{#1}_+}\xspace}
\newcommand{\Rs}{\ensuremath{R_\s}\xspace}
\newcommand{\As}{\ensuremath{A_{K,\s}}\xspace}

\newcommand{\T}[1][{}]{\ensuremath{\mathscr{T}_{\scriptscriptstyle #1}}\xspace}

\newcommand{\SjMain}[3]{\ensuremath{\mathscr{S}^{#1}_{#2}#3}\xspace}
\newcommand{\SjO}[2][\ut]{\SjMain{#2}{#1}{(\Omega)}}
\newcommand{\Sj}[2][]{\SjMain{#2}{\ut}{\brckt{#1}}}
\newcommand{\Sjn}[2][\ut]{\SjMain{#2}{#1}{}}

\newcommand{\So}[1][{}]{\ensuremath{\mathring{\mathscr{S}}_{\ut}^{#1}(\Omega)}\xspace}

\newcommand{\HVMain}[3]{\ensuremath{\tilde{\mathscr{S}}_{#1}^{#2}\brckt{#3}}\xspace}
\newcommand{\HV}[2][\ut]{\HVMain{#1}{#2}{\Omega}}

\newcommand{\HjMain}[3]{\ensuremath{\mathscr{H}^{#1}_{#2}#3}\xspace}
\newcommand{\Hj}[2][]{\HjMain{#2}{\ut}{\brckt{#1}}}
\newcommand{\HjO}[1]{ \HjMain{#1}{\ut}{(\Omega)} }

\newcommand{\hKR}[2]{\ensuremath{\mathsf{H}_{\!\scriptscriptstyle K\!R}^{#1}(#2)}\xspace}

\newcommand{\Kr}[1][0,q]{\ensuremath{\mathcal{K}^{#1}}\xspace}

\newcommand{\Acomma}[1]{\ifthenelse{\equal{#1}{}}{}{,#1}}

{\theorembodyfont{\rmfamily} }
{\theorembodyfont{\rmfamily} }
{\theorembodyfont{\rmfamily} }
{\theorembodyfont{\rmfamily} \theoremstyle{plain}  }

\newtheorem{thms}{Theorem}

\newcommand{\pf}{\noindent \textbf{Proof: }}
\newcommand{\epf}{\tiny \ensuremath{\hfill \blacksquare } \normalsize}
\newcommand{\eexm}{\tiny \ensuremath{\hfill \square } \normalsize}
\newcommand{\espf}{\small \ensuremath{\hfill \blacktriangledown } \normalsize}

\newcommand{\sName}{none}
\newcommand{\tName}{none}
\newcommand{\eeName}{none}
\newcommand{\eName}{none}

\newcounter{claim}
\newcommand{\setS}[1]{\label{S:#1}\renewcommand{\sName}{#1}}

\newcommand{\bgMain}[5][{}] {%
\renewcommand{\tName}{#4} 
\ifthenelse{\equal{#1}{}}{\begin{#2}\label{#3:#5:\tName}}
{\begin{#2}[#1]\label{#3:#5:\tName}}
}%

\newcommand{\bgT}[2][{}]{\bgMain[#1]{thm}{T}{#2}{\sName}{}}
\newcommand{\enT}{\end{thm}}

\newcommand{\bgP}[2][{}]{\bgMain[#1]{prop}{P}{#2}{\sName}{}}
\newcommand{\enP}{\end{prop}}

\newcommand{\bgL}[2][{}]{\bgMain[#1]{lemma}{L}{#2}{\sName}}
\newcommand{\enL}{\end{lemma}}

\newcommand{\bgE}[2][{}]{\bgMain[#1]{equation}{E}{#2}{\sName}{}}
\newcommand{\enE}{\end{equation}}

\newcommand{\bgD}[2][{}]{\bgMain[#1]{defn}{D}{#2}{\sName}{}}
\newcommand{\enD}{\end{defn}}

\newcommand{\bgC}[2][{}]{\bgMain[#1]{cor}{C}{#2}{\sName}{}}
\newcommand{\enC}{\end{cor}}

\newcommand{\bgR}[2][{}]{\bgMain[#1]{rem}{R}{#2}{\sName}{}}
\newcommand{\enR}{\end{rem}}

\newcommand{\enA}{\end{appl}}

\newcommand{\bgX}[2][{}]{\bgMain[#1]{exm}{X}{#2}{\sName}{}}
\newcommand{\enX}{\eexm \end{exm}}

\newcommand{\enQ}{\end{quest}}

\newcommand{\resetclaims}{\setcounter{claim}{1}}

\newcommand{\rfT}[2][\sName]{Theorem \ref{T:#1:#2}}
\newcommand{\rfC}[2][\sName]{Corollary \ref{C:#1:#2}}
\newcommand{\rfL}[2][\sName]{Lemma \ref{L:#1:#2}}
\newcommand{\rfP}[2][\sName]{Proposition \ref{P:#1:#2}}

\newcommand{\rfD}[2][\sName]{Definition \ref{D:#1:#2}}
\newcommand{\rfS}[1]{Section \ref{S:#1}}
\newcommand{\rfX}[2][\sName]{Example \ref{X:#1:#2}}

\newcommand{\rfE}[2][\sName]{\eqref{E:#1:#2}}

\newcommand{\bgEn}[2][\tName]{%
\renewcommand{\eName}{#1}
\renewcommand{\theenumi}{#2{enumi}}%
\begin{enumerate}}

\newcommand{\bgEnn}[2][\tName]{%
\renewcommand{\eeName}{#1}
\renewcommand{\theenumii}{$#2{enumii}$}%
\begin{enumerate}}

\newcommand{\bgEnS}[3][\tName]{\bgEn[#1]{#2} \addtocounter{enumi}{#3}}

\newcommand{\enEn}{\end{enumerate}}
\newcommand{\condition}[1]{#1\arabic }

\newcommand{\etem}{\item \label{i:\sName:\eName:\alph{enumi}}}
\newcommand{\eetem}{\item \label{ii:\sName:\eeName:\alph{enumii}}}

\newcommand{\rfi}[3][\sName]{(\ref{i:#1:#2:#3})}
\newcommand{\rfI}[1]{\rfi[\sName]{\tName}{#1}}
\newcommand{\rfii}[3][\sName]{(\ref{ii:#1:#2:#3})}

\newcommand{\rfCi}[3][\sName]{\rfC[#1]{#2} \rfi[#1]{#2}{#3}}

\begin{document}
\begin{center}
\large Boundary Regularity for the $\db$-Neumann Problem, Part 1 \normalsize\\

\hfill

Robert K. Hladky
\end{center}

\begin{abstract}
We establish sharp regularity and Fredholm theorems for the $\db$-Neumann problem on domains satisfying some non-generic geometric conditions. We use these domains to construct explicit examples of bad behaviour of the Kohn Laplacian: it is not always hypoelliptic up to the boundary, its partial inverse is not compact and it is not globally subelliptic. \end{abstract}

\begin{center}
\tableofcontents
\end{center}

\section{Introduction}\setS{ID}

In this paper we shall explore the $L^2$-existence and boundary regularity of solutions to the $\db$-Neumann problem on strictly pseudoconvex CR manifolds. Rather than attempt to solve the problem in general, we shall focus on a geometrically simple model class of examples for which positive results can be found.  We shall require that our domains possess a defining function $\rho$ depending upon the real and imaginary parts of a particular CR function $w$, and further that the level sets of $w$ produce a uniform foliation by compact, normal CR manifolds.

The first part of this condition is typical for most discussions of solvability for either the Kohn Laplacian $\boxb$ or the $\db$-complex on domains in CR manifolds. The study of the $\db$-Neumann problem initiated with Kuranishi  (\cite{Kuranishi1}, \cite{Kuranishi2} and \cite{Kuranishi3}), who established existence for a weighted Neumann problem on small balls, as part of his study of the embeddability of strictly pseudoconvex CR structures. More recently, Shaw has established unweighted $L^2$-existence results for small sets of CR manifolds embedded in $\cn{n}$ whose defining function satisfies both the condition above and a convexity constraint , see \cite{Shaw:B} or \cite{Shaw:P}. With the additional simplifying condition that the boundary has no characteristic points, Diaz has refined the techniques first employed by Kuranishi. In \cite{Diaz} he established that under certain curvature conditions, $L^2$ solutions exist with exact Sobolev regularity and estimates for a problem closely related to the $\db$-Neumann problem. His solutions are only guaranteed to meet the second Neumann boundary condition. Exact regularity refers to estimates of the type $\norm{\varphi}{k} \lesssim \norm{\boxb \varphi}{k}$. Diaz was interested in the tangential Cauchy-Riemann equations and his results are sufficient to show the existence of smooth solutions. However in general, the solutions exhibit a loss of Sobolev regularity. The author has not found any sharp estimates at the boundary for either the $\db$-Neumann problem or the tangential Cauchy-Riemann equations in the literature. On compact manifolds however, the Kohn Laplacian is well understood. In \cite{Folland:H} Folland and Stein introduced a new class of function spaces $\Sjn[]{k}$ and proved sharp estimates for the Kohn Laplacian in terms of these.

The analysis of the $\db$-Neumann problem is intricate as the operator \boxb is only subelliptic rather than elliptic. In addition the boundary conditions for the Neumann problem are non-coercive in the sense that the interior subelliptic estimates do not extend to the boundary of the domain $\Omega$. The presence of characteristic points on the boundary also complicates $L^2$ arguments enormously; the dimension of the horizontal space tangent to the boundary jumps. Both Kuranishi's and Diaz's argument for regularity involved the use of a subelliptic gain in directions tangent to the foliation by level sets of $w$. We shall work instead with the refined Folland-Stein spaces, but the underlying principle will remain the same. In a manner similar to these previous authors, we shall decompose the operator \boxb into pieces tangential and transverse to the foliation. However we shall work with global estimates on the compact leaves of a foliation and use local elliptic estimates in the transverse directions. 

The conditions we shall impose upon the CR manifolds and the subdomains are highly non-generic. This will yield many geometric simplifications, such as being able to choose the associated characteristic vector field to be tangential to the foliation, which will greatly facilitate the process of decomposing the operator \boxb. The boundary conditions split naturally and we will be able to make very specific identifications of portions of the operator with the Kohn Laplacians associated to the foliating manifolds.  Since each leaf of the foliation is compact, we can employ known results to obtain very precise estimates in directions tangential to the foliation.  In \cite{Tanaka} Tanaka constructed an explicit eigenvalue decomposition of \boxb on normal CR manifolds. We shall use this to decompose \boxb on the domain into an infinite family of elliptic operators on the hyperbolic plane. By studying these in detail and paying particular attention to uniformity of estimates over small hyperbolic balls, we shall be able to establish existence and sharp estimates for solutions to the $\db$-Neumann problem. 

Since our conditions are very restrictive, the negative implications of our results are perhaps of more general interest. For example although we show that \boxb has closed range as an operator on $L^2(D)$, its partial inverse $(1+\boxb)^{-1}$ is not a compact operator. This immediately implies that the $\db$-Neumann problem does not have global subelliptic estimates. A counter-example to compactness appeared in \cite{Shaw:P}, but the lower bound in the second estimate on page 162 appears to be incorrect. We shall provide a counter-example of our own in \rfS{BB}. In addition, the solutions to the \db-Neumann problem do not exhibit  a gain of Folland-Stein differentiability. We shall even show that hypoellipticity of $\boxb$  depends on Kohn-Rossi cohomology of the foliating manifolds. Since our examples  in many ways have the simplest possible geometry, it seems reasonable to conjecture that compactness and regularity gain fail generically. 

Despite the lack of global subelliptic estimates and the counter-examples to a gain of Folland-Stein regularity at the boundary, we are able to construct a family of spaces for which $1+\boxb$ is an isomorphism. The main theorem of Part 1 is

\begin{thms}\label{A}
 Let $\Omega$ be a smoothly bounded domain in a strictly pseudoconvex CR manifold of dimension $2n+1$ with $n \geq 3$ that satisfies the conditions of \rfS{MC}. 

\hfill

\noindent Suppose $1 \leq q \leq n-2$. Then on $(0,q)$-forms the operator
\[ 1+\boxb: \HV[]{k,2} \longrightarrow \SjO[]{k} \] is an isomorphism.
\end{thms}
When the Kohn-Rossi cohomology of the compact leaves vanishes at the $(0,q)$ level then we can replace $1+\boxb$ with $\boxb$.

The precise definition of the spaces and norms used here is given in Section \ref{S:BB}. We mention here that this theorem encodes exact regularity of solutions in the Folland-Stein spaces in all directions. Furthermore we obtain a full gain of two Folland-Stein derivatives for all directions in the interior and in directions tangent to the foliation at the boundary. In particular, hypoellipticity at the boundary of solutions for $1+\boxb$ is implied.

The method used to construct our class of examples unfortunately precludes the possibility of our domains possessing characteristic points. However in the special case of the Heisenberg group and a foliation by spheres, the arguments can be extend to cover characteristic domains. This will be explored in detail in Part 2 for the unit ball of the Heisenberg group.

An important application of the \db-Neumann problem is to solving the inhomogeneous tangential Cauchy-Riemann equation.  Our main theorem yields an existence and regularity theory for this problem. Our result is as follows:

\begin{thms}
Under the same conditions as Theorem \ref{A}, for any $(0,q)$-form $\vs \in \lO[] $ such that $\db \vs = 0$ and $\vs \bot \Ker{\boxb}$ the equation
\[ \db \varphi = \vs \] is uniquely solvable for $\varphi \in \HV[]{0,1}$. Furthermore, if $\vs \in \SjO[]{k}$ then $\varphi \in \HV[]{k,1}$ and there is a uniform estimate 
\[ \norm{\varphi}{\HV[]{k,1}} \lesssim \norm{\vs}{\SjO[]{k}}.\]
\end{thms}

Again the precise definitions of the spaces involved is put off until \rfS{BB}. However we note that this encodes exact regularity in the weighted Folland-Stein spaces with a slight  gain in directions tangential to the foliation. This is sufficient to establish solutions globally smooth up the boundary when $\vs$ is itself smooth.

Both Parts of this paper based on the Author's Ph.D. thesis at the University of Washington and benefited enormously from the guidance and criticism given by John M. Lee.

\section{Pseudohermitian Manifolds and the Kohn Laplacian}\setS{PS}

\bgD{CR}
An abstract (hypersurface-type) CR-structure on a (2n+1)-manifold $M$ is a complex $n$-dimensional subbundle  $\crT \subset \cn{}TM$ such that 
\bgEn{\roman}
\etem $\crT \cap \acrT =\{0\}$ where $\acrT = \overline{\crT}$,
\etem $[C^\infty(\crT),C^\infty(\crT)] \subset C^\infty(\crT)$. 
\enEn
\enD
An odd dimensional manifold $M$ paired with a CR-structure $\crT$ is referred to as a CR-manifold.

The real expression of a CR-structure consists of the real space $H:= \text{Re }(\crT \oplus \acrT) \subset TM$ together with a real automorphism $J:H \to H$ such that $\crT$ is the $+i$-eigenspace and $\acrT$ the $-i$-eigenspace for $J$. The integrability condition is then equivalent to the vanishing of the Niunhuis tensor
\[ N_J(X,Y):= [X,JY]+[JX,Y]- J \left( [X,Y]-[JX,JY] \right).\]
The real expression of the CR-structure on $M$ is a real subbundle $H$ of $TM$ with dimension $2n$. Thus $H^{\perp} =\{\ua \in T^*M: \ua(X)=0 \text{ for all } X\in H\}$ is a 1-dimensional line bundle on $M$. Suppose that $H^{\perp}$ is trivial and admits a non-vanishing global section $\theta$. The choice of the 1-form $\theta$ is not canonical for the CR-structure. However, since $\theta$ is determined up to multiplication by a non-vanishing $C^\infty$ function, much can be derived from $\theta$ that is actually a property of the underlying CR-structure. 

Alternatively, a triple $(M,J,\theta)$ is referred to as a \textit{pseudohermitian manifold} if the $1$-form $\theta$ is non-vanishing and the endomorphism $J:\Ker{\theta} \to \Ker{\theta}$ satisfies $J^2= -\one$ and $N_J =0$. Set $H=\Ker{\theta}$ and $\crT = \{ X-iJX: X \in H\}$ then $(M,\crT)$ is the \textit{underlying CR structure} for the pseudohermitian manifold.

\bgD{Levi}
The Levi form for a pseudohermitian manifold $(M,J,\theta)$ is defined by
\[ L_\theta(X,Y) =d\theta(X,JY) \quad \text{for $X,Y \in C^\infty (H) $}.\]
\enD 

While $L_\theta$ depends on the choice of contact form $\theta$, non-degeneracy and number of non-zero eigenvalues depend only on the underlying CR-structure. In particular, we will say that a pseudohermitian (or CR) manifold is \textit{strictly pseudoconvex}  if the eigenvalues for its associated Levi form are either all positive or all negative. By multiplying $\ut$ by $-1$ if necessary we can always assume that they are all positive.

It is well known that if the pseudohermitian manifold $(M,J,\theta)$ has a non-degenerate Levi form then there exists a unique (real) vector field $T$ such that $\theta(T)=1$ and $T \lrcorner d\theta =0$. 
This vector field, $T$, is referred to as the \textit{characteristic field} for the pseudohermitian form $\theta$. We obtain a geometric interpretation of the Levi form. Specifically, since $\theta$ annihilates $H$, we see that $L_\theta(X,Y)$ is the component of $-[X,JY]$ in the characteristic direction. 

The characteristic field allows us to extend $J$ to a linear bundle map $TM \to TM$ by requiring $JT=0$. A strictly pseudoconvex pseudohermitian structure thus induces a Hermitian metric on  $M$ defined by 
\[ h_\ut(X,Y) = d\theta(X,J\bt{Y}) + \theta(X) \theta(\bt{Y}).\]
In particular $h_\ut(X,Y) =L_\theta(X,Y)$ for $X,Y \in H$. When restricted to $TM$ this defines a Riemannian metric. We shall call this the Webster metric. Suppose $Z,W \in \crT$, then \[h_\ut(Z,\bt{W})=d\theta(Z,iW)= -i\theta([Z,W])=0\] by the integrability condition. Thus the Webster metric is compatible with the CR structure in the sense that $\crT$ and $\acrT$ are orthogonal. Indeed there is an orthogonal decomposition
\[ \cn{}TM = \crT \oplus \acrT \oplus \cn{}T.\]  
For a smooth vector field $X$ we shall denote by $X\upp$ and $X\dupp$ the projections of X to \crT and \acrT respectively. 
A strictly pseudoconvex pseudohermitian structure possesses a canonical connection. This allows us to intrinsically define a variety of Sobolev type spaces in addition to providing a useful computational tool.

\bgL{Connection}
If $(M,J,\theta)$ is strictly pseudoconvex, there is a unique connection $\nabla$ on $(M,J,\theta)$ that is compatible with the pseudohermitian structure in the sense that $H$, $T$, $J$ and $d\theta$ are all parallel and the torsion satisfies
\begin{align*}
\text{Tor}(X,Y)&=d\theta(X,Y)T\\
\text{Tor}(T,JX)&= -J\text{Tor}(T,X)
\end{align*}
\enL

This formulation of the connection was developed by Tanaka \cite{Tanaka}. An alternative formulation in terms of a coframe was independently derived by Webster  \cite{Webster}.  Tanaka's proof is constructive and is based upon the useful formula
\bgE{Connection}
\nabla_{X\dupp} Y\upp = [X\dupp,Y\upp]\upp.
\enE
The curvature $R^\nabla$ of the connection is the $(3,1)$-tensor defined by
\bgE{Curvature}
R^\nabla \varphi(X,Y) = \nabla_X \nabla_Y \varphi - \nabla_Y \nabla_X \varphi - \nabla_{[X,Y]} \varphi.
\enE

The Webster metric associated to $\ut$ induces a metric on every tensor bundle over $M$. The Webster-Tanaka connection also induces a connection on every bundle. We shall denote by $\mathcal{T}^kM$ the space of contravariant $m$-tensors over $M$. For all $k$ the Webster-Tanaka connection then induces a natural map
\[ \nabla: C^\infty(\mathcal{T}^m \otimes \Lambda^{k}M) \to C^\infty(\mathcal{T}^{m+1} \otimes \Lambda^{k}M).\]

The pseudohermitian form $\ut$ induces a canonical  volume form for the triple $(M,J,\theta)$, given explicitly by $dV= \theta \wedge (d\theta)^n$. This yields a canonical  global $L^2$ norm and inner product.  We can now intrinsically describe a family of Sobolev type spaces for the pseudohermitian manifold $(M,J,\ut)$. For some choice of non-negative integers $k$ and $m$, let $E$ be the bundle $\mathcal{T}^m M \otimes \Lambda^{k} M$. The $L_\ut^2$-Sobolev spaces for such bundles $E$ are then inductively defined by
\bgE{Sobolev}
\begin{split}
\Hj[E]{0} &= L^2_\ut(E)\\
\Hj[E]{j} &= \left\{\varphi \in \Hj[E]{j-1}: \nabla \varphi \in \Hj[\mathcal{T}^1M \otimes E]{j-1} \right\}, \qquad \text{$(j>0)$} 
\end{split}
\enE
with norms defined by
\[ \norm[2]{\varphi}{\Hj[E]{j})} = \sum\limits_{k \leq j} \norm[2]{ \nabla^k \varphi}{L_\ut^2(\mathcal{T}^k M \otimes E)}. \] 
If $m=0$ we shall denote these spaces simply by $\Hj[M]{j}$.

These spaces do not provide optimal results for the analysis of the Kohn Laplacian as it is not fully elliptic; it's only first order in the characteristic direction. The Folland-Stein spaces were introduced in \cite{Folland:H} to provide more refined regularity results. They are defined intrinsically on $(M,J,\ut)$ as follows. For a differential form $u$, we can decompose \[\nabla \varphi = \nabla \upp \varphi + \nabla \dupp \varphi + \nabla_T  \varphi \otimes \theta.\]   Here we define  $\nabla\upp$ by $\nabla\upp \varphi (\cdot, X)=\nabla \varphi (\cdot, X\upp)$ and $\nabla\dupp$ by $\nabla\dupp \varphi(\cdot, X)=\nabla \varphi(\cdot,X\dupp)$.
Now set 
\[\NabH \varphi := \nabla \upp \varphi + \nabla\dupp \varphi= \nabla \varphi - \nabla_{T} \varphi \otimes \theta.\]
The Folland-Stein spaces are now defined in an analogous fashion to the Sobolev spaces using the restricted operator $\NabH$ instead of $\nabla$. More concretely
\bgE{FollandStein}
\begin{split}
\Sj[M]{0} &= L^2_\ut(M)\\
\Sj[M]{j} &= \left\{ \varphi \in \Sj[M]{j-1}: \NabH \varphi \in \Sj[M]{j-1} \right\}, \qquad \text{for $j>0$} 
\end{split}
\enE
with norms
 \[\norm[2]{\varphi}{\Sj[M]{j}} = \sum\limits_{k \leq j} \norm[2]{\left( \NabH \right)^k \varphi}{L_\ut^2(M)}. \]
A lot of notation has been suppressed in this definition. A more careful description would follow the lines of \rfE{Sobolev}.

The spaces $\HjO{j}$ 
and $\SjO{j}$ 
can be defined for open subsets $\Omega \subset M$ simply by replacing $M$ with $\Omega$ in the definitions. 

When working locally in some well-behaved frame, the above definitions can be replaced by more concrete descriptions involving derivatives of functions. The nature of the argument expounded in Part 2 of this paper requires careful analysis on unbounded domains, where the most convenient frames are singular or degenerate. The intrinsic definitions given above allow provide a solid foundation for dealing with this situation. 

Usually when we refer to the Folland-Stein spaces we shall impose an additional constraint on the forms under consideration. This constraint is the notion of bidegree.

For the remainder of this section we suppose $(M,J,\theta)$ is a strictly pseudoconvex pseudohermitian manifold. Set $\Lambda_\ut^{0,1}M=\{ \varphi \in \cn{}T^*M: \varphi=0 \text{ on } \crT \oplus \cn{}T\}$ and let $\Lambda_\ut^{1,0}M$ be the orthogonal complement to $\Lambda_\ut^{0,1}M$ in $\cn{}T^*M$. It should be stressed that these definitions are asymmetric. We extend to higher degree forms by setting $\Lambda_\ut^{p,q}M = \Lambda^p\left( \Lambda_\ut^{1,0} M \right)\otimes \Lambda^q \left(\Lambda_\ut^{0,1}M\right)$. The space of degree $k$ complex covector fields on $M$ then admits the following orthogonal decomposition \[ \cn{}\Lambda^{k} M = \bigoplus\limits_{p+q=k} \Lambda_\ut^{p,q}M.\] Denote the orthogonal projection $\cn{}\Lambda^{p+q}M \to \Lambda_\ut^{p,q}M$ by $\pi_\ut^{p,q}$ and define
\bgE{concrete}
\db = \pi_\ut^{p,q+1} \circ d.
\enE
Then $\db$ maps $C^\infty(\Lambda_\ut^{p,q}) \to C^\infty(\Lambda_\ut^{p,q+1})$. It should be remarked that this definition depends upon the pseudohermitian structure and is not canonical for the underlying CR structure. Using the language of holomorphic vector bundles and quotients, it is possible to construct an operator depending solely on the CR structure that reduces to our definition once a pseudohermitian form is chosen. However, for the purposes of this paper the concrete version offered here will suffice. It can also be computed in terms of the Webster-Tanaka connection.

\bgL{Computedb}
For any smooth $(0,q)$-form $\varphi$
 \[\db \varphi = (-1)^{q}(q+1)\textup{Alt } (\nabla\dupp \varphi).\]
\enL

\pf
We shall prove the lemma for $(0,1)$-forms. The general case then follows from the observation that both operators obey the same product rules. For a $(0,1)$-form $\varphi$
\begin{align*}
2\text{Alt} (\nabla \dupp \varphi) (X,Y)&=\nabla \varphi (X,Y\dupp)- \nabla \varphi (Y,X\dupp)= \nabla_{Y\dupp}\varphi ( X) - \nabla_{X\dupp} \varphi(Y)\\
&= Y\dupp \varphi(X) - X\dupp \varphi(Y)  - \varphi(\nabla_{Y\dupp} X - \nabla_{X\dupp } Y)\\
&= Y\dupp \varphi(X\dupp) -X\dupp\varphi(Y\dupp)- \varphi(\nabla_{Y\dupp} X\dupp - \nabla_{X\dupp } Y\dupp)\\
&= Y\dupp \varphi(X\dupp) -X\dupp\varphi(Y\dupp) - \varphi([Y\dupp,X\dupp] +Tor(Y\dupp,X\dupp))\\
&= -\db \varphi (X,Y).
\end{align*}
Here we have made implicit use of the fact that $\crT$ and $\acrT$ are parallel and that the torsion of two sections of \acrT vanishes.

\epf

It follows immediately from the definitions that $\db \circ \db =0$. Thus $\db$ defines a complex of differential forms on $M$. The associated cohomology is known as the Kohn-Rossi cohomology and is denoted by $\hKR{p,q}{M}$. A key tool for studying these groups is the Kohn Laplacian.
\bgD{Kohn}
If $(M,J,\theta)$ is a strictly pseudoconvex structure then the formal Kohn Laplacian is defined by
\[ \boxb = \db \vtb + \vtb \db \]
where \vtb denotes the formal adjoint of \db with respect to the canonical $L^2$ inner product on $(M,J,\ut)$.
\enD

\bgT{compact} Let $(M,J,\ut)$ be a compact, strictly pseudoconvex pseudohermitian manifold of dimension $2n-1$. 
\bgEn{\alph}
\etem If $1 \leq q \leq n-2$ then $\boxb$ is a self-adjoint, Fredholm operator on $L^2(\Lambda_\ut^{p,q}M)$ and there is an orthogonal decomposition
\[ \begin{split}L^2(\Lambda_\ut^{p,q}M) &= \rng{\boxb} \oplus \Ker{\boxb}\\
&= \rng{\db} \oplus \rng{\vtb} \oplus \Ker{\boxb}.\end{split}\]
The operator $\boxb$ is subelliptic. Therefore $\Ker{\boxb}$ is finite dimensional  and $(1+\boxb)^{-1}$ is a compact, bounded operator on $L^2$. The cohomology group $\hKR{p,q}{M} \cong \Ker{\boxb}$. Furthermore the operator $1+\boxb$ is an isomorphism from \Sj[\Lambda^{p,q}_\ut M]{k+2} to \Sj[\Lambda^{p,q}_\ut M]{k} for all $k \geq 0$. 
\etem If $q=0$ then $\boxb$ is self-adjoint and has closed range as an operator on $L^2(\Lambda_\ut^{p,0}M)$ and there is an orthogonal decomposition
\[ L^2(\Lambda_\ut^{p,0}M) = \rng{\boxb} \oplus \Ker{\boxb}.\]

\enEn
\enT

\bgR{dimension}
The reader is cautioned that the dimension assumed in the theorem is $2n-1$ rather than the $2n+1$ of \rfD{CR}. When we are adding the supposition that the manifold is compact we shall always adopt this drop of dimension, whereas if we are not presupposing compactness we shall continue to use $2n+1$. This is to ensure compatibility of results when we are working with a foliation of a domain by compact CR manifolds of codimension $2$.
\enR

\renewcommand{\tName}{compact}
\pf
The results of \rfI{a} are the culmination of the work by Folland and Stein in \cite{Folland:H}. The self-adjointness of \rfI{b} is similar to that of \rfI{a}. To show closed range we first note that on $(p,0)$-forms $\boxb = \dbs \db$ and $\Ker{\boxb} =\Ker{\db}$. Clearly $\rng{\boxb} \subset \Ker{\db}^\bot $. Suppose $\varphi \bot \Ker{\db}$. By \rfI{a} we can write $\db \varphi = \boxb \psi + \varsigma$ where $\varsigma \in \Ker{\boxb}$ and $\psi \bot \Ker{\boxb}$.  Using the orthogonality properties it follows that $\varsigma =0$ and $\varphi = \dbs \psi$. Now we again use \rfI{a} to write $\psi = \boxb \tau$. Then $\varphi = \dbs \boxb \tau = \dbs \db \dbs \tau = \boxb \dbs \tau$. Thus $\rng{\boxb} = \Ker{\db}^\bot$.

\epf

\bgR{OtherP}
Throughout the literature, most computations and arguments concerning $\db$ are conducted under the assumption that $p=0$. When the bundle $\Lambda^{1,0}_\ut M$ is holomorphically trivial, it is easy to pass to the general case. Since $\cn{}TM/\acrT $ is a holomorphic bundle, local results always simply to the case $p>0$.  If $M^{2n-1}$  is globally CR embeddable as a hypersurface in $\cn{n}$ then the embedding functions  provide a global trivialisation. \enR

For a compact manifold, the formal and actual $L^2$ adjoints of \db are equal. However for manifolds with boundaries the issue of boundary conditions arises and we must make a subtly different definition to recover self-adjointness for the operator.

Suppose $\Omega$ is a bounded open set in $M$ with smooth boundary. We can restrict our complexes to forms defined on $\Omega$. We extend $\db$ to its maximal $L^2$  closure, also denoted \db, and define $\dbs$ to be the $L^2$-adjoint of this extended operator. We then define the Kohn Laplacian for $\Omega$ by
\[\dom{\boxb} = \{\varphi \in \dom{\db} \cap \dom{\dbs}: \text{ $\db \varphi \in \dom{\dbs}$ and $\dbs \varphi \in \dom{\db}$} \}\]
and for  $\varphi \in \dom{\boxb}$
\[ \boxb \varphi = \db  \dbs \varphi + \dbs  \db \varphi.\] 
This operator is  self-adjoint as an operator on $L^2(\Lambda_\ut^{p,q}\Omega)$ (see \cite{Shaw:B}). For forms contained in $\dom{\boxb}$ the operator agrees with the formal version defined above. The $\db$-Neumann problem on $\Omega$  is then to decide when the equation $\boxb u =f$ on $\Omega$ can be solved for $u \in \dom{\boxb}$  and obtain optimal regularity results. It is worth emphasising that there are boundary constraints on any solution $u$ for it to lie in $\dom{\boxb}$.  From the view point of CR geometry as opposed to pseudohermitian geometry we would also be free to choose an appropriate $\ut$. 

The analysis of this problem is difficult for several reasons. The operator is not elliptic as it has only limited control over the characteristic direction. However, the Folland-Stein spaces were constructed to address precisely this.  The characteristic vector field $T$ can be written as a commutation of vector fields from $\acrT$. Thus $T$ is second order as an operator in the Folland-Stein setting. Although \boxb is only subelliptic (see \cite{Shaw:B}), it is fully elliptic in the Folland-Stein directions. 
A second problem  is that the boundary conditions are non-coercive in a sense to be made more precise later. There is also a third problem related to the geometry of the boundary \dO.

\bgD{CharPoint}
A point $x \in \dO$ is a characteristic point for $\Omega$ if the boundary is tangent to the distribution $H$ at $x$, i.e. $T_x\dO =H_x$.
\enD

At all non-characteristic points the tangent space to $\dO$ intersects $H$ transversely with codimension $1$. Thus at characteristic points there is a jump in the dimension of this intersection. This phenomenon makes obtaining $L^2$ estimates difficult near these points. 

All positive results for this problem have required strict conditions on the geometry of the boundary of $\Omega$ and have returned non-sharp boundary regularity. See for example the work by Ricardo Diaz \cite{Diaz} or Mei-Chi Shaw \cite{Shaw:P}. In particular very little is known about regularity when the domain possesses characteristic points.

\section{Compact, Normal Pseudohermitian Manifolds}\setS{NM}

Under a certain geometric condition on a compact,  pseudohermitian manifold it is possible to go much further with the study of the Kohn Laplacian. If some key operators commute then the techniques of harmonic analysis can be employed to construct a detailed and useful eigenform decomposition. The reader is cautioned that in this section we shall assume dimension $2n-1$ for our manifolds rather than $2n+1$ as earlier. The reason for this shall soon become apparent. This section is a summary of work done by Tanaka \cite{Tanaka} organised in a form more convenient for our current purposes.

\bgD{Normal}
A strictly pseudoconvex pseudohermitian manifold $(N,J,\ut)$ with characteristic field $T$ is said to be normal if either of the following equivalent conditions hold
\bgEn{\condition{N}}
\etem $T$ is analytic, i.e. $[T,C^\infty(\crT)] \subseteq C^\infty(\crT)$.
\etem $\mathcal{L}_T J=0$, i.e. $[T,JX]=J[T,X]$ for all $X \in C^\infty(H)$.
\enEn
\enD
For a more thorough study of normal strictly pseudoconvex manifolds the reader is directed to \cite{Tanaka}. Here we shall immediately focus in on the specific parts of the harmonic analysis that we shall need later.

\bgD{Lambda}
The differential operator $\T$ is defined on smooth differential forms
\[ \T \varphi = -i \nabla_T \varphi\]
and extended to be a maximal, closed unbounded operator on each $L^2(\Lambda_\ut^{p,q}N)$.
\enD
Note that the normality condition implies that $\T$ preserves bidegree so this definition makes sense.

\bgL{Basic}
If $(N,J,\ut)$ is a normal, strictly pseudoconvex pseudohermitian manifold then
\bgEn{\alph}
\etem $\nabla_T X = [T,X]$
\etem $[\T,\db]=0$.
\etem $T \lrcorner R^\nabla =0$
\enEn
\enL

\pf Since $N$ is normal if $X \in C^\infty(H)$ then
\[ \text{Tor}(T,JX)= \nabla_T JX - [T,JX] = J\nabla_T X - J[T,X] = J \text{Tor}(T,X).\]
But by the defining properties of the Webster-Tanaka connection we have $\text{Tor}(T,JX)=-J\text{Tor}(T,X) $. Thus $\text{Tor}(T,JX)=0$ for all such $X$. Thus $T \lrcorner \text{Tor} =0$. Then for any vector field $X$
\[ \nabla_T X = \text{Tor}(T,X) + [T,X] = [T,X].\] 

To see \rfI{b} we note that \rfI{a} and \rfi{Normal}{a} imply that $\T$ preserves bidegree. Thus for a $(p,q)$-form $\varphi$, $[\T,\db]\varphi $ is the projection onto $(p,q+1)$-forms of $[\T,d]$. But $\T = -i \mathcal{L}_T$ and all Lie derivatives commute with the exterior derivative.

For \rfI{c} we note that the tensors $\ut$, $d\ut$ and $J$ are all invariant under $T$ and so the Webster-Tanaka connection is also. Thus
\[ \mathcal{L}_T \nabla_X \varphi = \nabla_X \mathcal{L}_T \varphi + \nabla_{[T,X]} \varphi.\]
The result then follows from \rfI{a}.

\epf

The commutation property exhibited in \rfI{b} allows for the development of a very useful harmonic theory on compact, normal manifolds. 

\bgT[Tanaka]{qNonZero}
Suppose $N$ is a compact, normal, strictly pseudoconvex pseudohermitian manifold of dimension $2n-1$ with $n \geq 3$. Then for $0 \leq q \leq n-2$ there is a countable collection of smooth $(0,q)$ forms $\V = \bigcup\limits_{\Gamma,\ul \in \rn{}}  \V(\Gm[],\ul)$ such that
\bgEn{\alph}
\etem \V is an orthonormal basis for $L_{(0,q)}^2(N)$.
\etem For $\s \in \V(\Gm[],\ul)$ we have $\boxb \s = \Gm[] \s$ and $\T \s = \ul \s$.
\etem For $\s \in \V$,  $\Gm \geq -(n-q-1)\uls$.
\etem There exists $\Go>0$ such that the space $K_1^q := \{\Gm: \s \in \V\}$ is a subset of $\{0\} \cup [\Go,\infty)$. 
\etem For each $\Gm[] \in K_1^q$ the image of $\V(\Gm[]) := \bigcup\limits_\ul \V(\Gm[], \ul)$ under the map $\s \mapsto \uls$ is discrete in $\rn{}$.
\etem The space $K^q_2 := \{(\Gm,\uls):\s \in \V\}$ is a discrete subspace of $[0,\infty) \times \rn{}$. 
\enEn
Furthermore if $q>0$ then 
\bgEnS{\alph}{6}
\etem $K_1^q$ is discrete.
\etem There exists a constant $C>0$ such that $\snorm{\uls}{} \leq C (1 + \Gm)$ for all $\s \in \V$.
\etem Each $\V(\Gm[],\ul)$ is a finite (possibly empty) set.
\enEn
\enT

\pf (sketch)

When $q>0$ most parts of this theorem are fairly standard consequences of the subelliptic theory for  the positive operator $\boxb$ on compact strictly pseudoconvex manifolds. Again the reader is referred to \cite{Tanaka} for details, although they should be cautioned that the sign convention adopted here for \T is the opposite to that assumed there. Part \rfI{c} follows from an identity also shown in \cite{Tanaka}, that
\bgE{Bar}
 \boxbt \varphi - (n-q-1) \T \varphi = \boxb \varphi
 \enE
for smooth $(0,q)$-forms ($0 \leq q \leq n-1$). Here  $\boxbt \varphi := \Bt{ \boxb \bt{\varphi}}$.  Part \rfi{qNonZero}{g} follows immediately from a result of Folland and Stein \cite{Folland:H} that $1+ \boxb$ is an isomorphism  from $\Sj[N]{2}$ to $L^2(N)$.

We note in that the existence of $\Go$ in \rfi{qNonZero}{d}, indicating the presence of a spectral gap above the kernel for \boxb, is implied by \boxb having closed range.

When $q=0$ we still have closed range, but since $\boxb$ is no longer subelliptic we shall instead base the argument around the self-adjoint operator $P = 1 + \boxb + \boxbt$.  Now it is easy to see that $P = 1+\triangle_b$ where $\triangle_b$ is the sub-Laplacian of \cite{Lee:Y}.  Thus $P$ is subelliptic on functions. Furthermore $P+\T^2$ is strongly elliptic and $P$ commutes with $\T$ when applied to smooth functions. We can therefore construct a smooth orthonormal basis \V[0] for $L^2_{(0,0)}(N)$ such that each $\s \in \V[0]$ is an eigenfunction for both $P$ and \T. In addition the pairs of eigenvalues $(p(\s),\uls)$ form a discrete subset of \rn{2}. However we can use \rfE{Bar} to see that $P = 1+2\boxb + (n-1)\T$. Thus each $\s \in \V[0]$ is  also an  eigenfunction for $\boxb$.  Furthermore the map $(x,y) \mapsto (x-1-(n-1)y,y)$ is an affine transformation of \rn{2} that taking pairs of eigenvalues $(p(\s),\uls)$ for $P$,$\T$ to pairs of eigenvalues $(\Gm,\uls)$ for $\boxb$, $\T$. Combined with the observation that \boxb is positive,  this establishes the discreteness result of \rfi{qNonZero}{f}. All the other properties now follow from standard arguments and \rfE{Bar}. 

\epf

\bgD{BarGamma}
For $\s \in \V[0]$ we define $\Bt{\Gamma}(\s)$ by
\[ \boxbt \s = \Bt{\Gamma}(\s) \s.\]
Then for $\s \in \V$  we set \[ G(\s) =\begin{cases} \big(1+\Gm +\Bt{\Gamma}(\s) \big)^{1/2}, \quad &\text{if $\s \in \V[0]$},\\ \big(1+\Gm \big)^{1/2}, \quad &\text{if $\s \in \V$, $q>0$.}\end{cases}\] We can then define an unbounded linear operator $G$ on $L^2$ which maps each $\s$ to $G(\s)\s$ and is extended by linearity. 
\enD
These definitions make sense because of \rfE{Bar}, but the reader is cautioned to note that $\Bt{\Gamma}(\s)$ is a real function and not the conjugate of $\Gm$.

\bgL{FollandStein}
Under the same conditions as \rfT{qNonZero} the linear operator  $G$ is an isomorphism from $\Sj[N]{k+1}$ to $\Sj[N]{k}$ for all $k \geq 0$.
\enL

\pf
The case $q>0$ follows from \rfT{qNonZero}. When $q=0$ the result follows from the fact that $1+\triangle_b$ acts as an isomorphism between the requisite Folland-Stein spaces. The reader is referred to \cite{Lee:Y} for details.

\epf

We conclude this section by noting some important examples of normal pseudohermitian manifolds.

\bgX[Unit Sphere]{Sphere}
Let $N=\sn{2n-1}$ be the unit sphere in \cn{n}. We impose a strictly pseudoconvex structure on $N$ by setting
\[ \ut =\frac{i}{2} \left( z^k d\bt{z}^k - \bt{z}^k dz^k \right)\]
and taking $J$ to be the map naturally induced from the holomorphic structure of \cn{n}. The characteristic vector field is then given by \[T = i\left({z}^k \frac{\partial}{\partial {z}^k} - \bt{z}^k \frac{\partial}{\partial \bt{z}^k} \right)\]
which is just the infinitesimal generator of the natural $U(1)$ action on $N$. The $U(1)$ action acts by contact diffeomorphisms so $N$ is normal. In \cite{Folland:S} Folland conducted a detailed study of the harmonic structure of \boxb. We shall exploit this heavily in Part 2 of this paper. Here we note that $\Lambda^{1,0}_\ut N$ is holomorphically trivial and the Kohn-Rossi cohomology for $N$ at the $(p,q)$ level vanishes when $1 \leq q \leq n-2$.
\enX

\bgX[Heisenberg Group]{Heisenberg}
This example is non-compact but will be of considerable importance in Part 2. Set $\hn{2n+1}=\{ (t,z)\in \rn{} \times \cn{n} \}$. A contact structure can be imposed by setting $\ut = dt + iz^k d\bt{z}^k -i\bt{z}^k dz^k$.  Define a function $w$ by $w:=t + iz^k \bt{z}^k$. The operator $J$ is then  induced from the embedding of \hn{2n+1} into \cn{n+1} given by $(t,z) \mapsto (w,z)$. The characteristic vector field is given by $T=\frac{\partial}{\partial t}$ which generates the translations along the $t$-axis.

A second normal structure with the same underlying CR structure can be imposed on the subset $\hn{2n+1} -\{z=0\}$ by defining $\tilde{\ut} = (2 z^k \bt{z}^k)^{-1} \ut$. Although singular at the set $\{z=0\}$ this pseudohermitian has the analytically pleasing properties that it is preserved by both dilations and translations. The characteristic field is given by $T=    i\left(z^k \frac{\partial}{\partial {z}^k} - \bt{z}^k \frac{\partial}{\partial \bt{z}^k} \right)$. We then note that the level sets of the CR function $w$ are all naturally isomorphic to the unit sphere as pseudohermitian manifolds.
\enX

\bgX{Kahler}
Let $V$ be a compact K\"ahler manifold and $E$ a holomorphic line bundle over $V$ that is negative in the sense that the first Chern class of $E$ can be represented by a negative $(1,1)$-form $\omega$. Recall that a $(1,1)$-form $\omega$ is negative if the symmetric form $\omega(X,JY)$ is negative definite. It follows from a theorem presented by Morrow and Kodaira (\cite{Morrow}, Theorem 7.4) that there is a connection on $E$ with curvature  $i\omega$. In the $U(1)$-principle bundle $M$ associated to $E$, this connection defines a choice of horizontal distribution.  The Lie algebra for $U(1)$ is naturally identifiable with $i \rn{}$. Thus $\ut$, defined to be $i$ times the connection $1$-form, is a non-vanishing real valued $1$-form on $M$. The projection from the horizontal distribution to the tangent space of $V$ induces a natural $J$ operator for $M$.  Now under suitable identifications $d\ut = - \omega$ on $H$. Thus the structure $(M,J,\ut)$ is strictly pseudoconvex. The characteristic field $T$ is easily seen to be the infinitesimal generator of $1$-parameter family given by the right action of $e^{-it}$. These extend to holomorphic transformations of $E$ and so we get $\mathcal{L}_T J=0$. Hence $M$ is normal. 
\enX
\section{The Class of Model Examples}\setS{MC}

In this section we introduce the class of model domains for our study of the $\db$-Neumann problem. We begin by assuming that $(N,J_N,\ut_N)$ is a compact, normal strictly pseudoconvex pseudohermitian manifold of dimension $2n-1$ with $n \geq 3$. We shall use the notation \[\hy{2}:= \{ w =t+is \in \cn{}: s>0\}\]  to denote the upper halfplane model of $2$ dimensional hyperbolic space.

\bgD{MTilde}
Set $M$ to be the smooth product manifold $\hy{2} \times N$. Let $\pi:M \to N$ be the natural projection onto $N$ and $w:M \to \hy{2}$ the natural projection onto \hy{2}. 
\enD

We shall also use $w$ as both a holomorphic coordinate on \hy{2} and as a complex valued function on $M$. We identify all covariant tensors on $N$ with their pullbacks via $\pi$. For any vector field $X$ on $N$ we denote by $\tilde{X}$ the unique vector field on $M$ such that $\pi_* X = N$ and $Xw=X\bt{w}=0$. 

\bgD{Y}\hfill

\bgEn{\roman}
\etem Define  $Y$ to be the unique (smooth) vector field on $M$ such that $\pi_* Y =-\frac{i}{2} T_N$ and $w_* Y = 2is \frac{\partial}{\partial w}$.
\etem A differential form $\varphi$ is said to be tangential if $Y \lrcorner \varphi = \bt{Y} \lrcorner \varphi =0$.
\etem Set $\ut^0 = \frac{1}{2is} dw$ and $\ut^\bt{0} = \Bt{\ut^0}$. We identify these with $w^* \ut^0$ and $\bt{w}^* \ut^\bt{0}$ on $M$.
\enEn
\enD

\bgL{Example}
Construct a non-vanshing real $1$-form $\ut$ on $M$ by
\[ \ut = \ut_N + \frac{i}{2} \ut^0 -\frac{i}{2} \ut^\bt{0}\] and define a distribution in $\cn{}M$ by
\[ \crT M =  \{ \tilde{X} : X \in \crT N \}  \oplus \text{span}_{\cn{}} \langle Y \rangle.\]
Then $(M,\crT M)$ is a strictly pseudoconvex CR manifold that admits $\ut$ as a pseudohermitian form. Furthermore the characteristic vector field $T$ is given by $T=\widetilde{T_N}$.
\enL

\pf
For convenience we shall denote the first component of $\crT M$ by $(\crT M)^\top$. The distribution $\crT M$ clearly has the correct dimension for the pair $(M,\crT M)$ to be CR. It just remains to show that its sections are closed under the Lie Bracket. This is true for each separate component of  $\crT M$ so it just remains to show this for $[\tilde{X},Y]$ for any section $X$ of $\crT N$. Now by the defining properties of normal pseudohermitian manifolds we see that
\[ \pi_* [\tilde{X},Y] = [ X , \pi_* Y] = [ X,-\frac{i}{2}T_N] \in \crT N.\] In addition $Y\bt{w}=0$ and $Yw  = 2i \imag{w} = w- \bt{w}$. Thus $[\tilde{X},Y]w = [\tilde{X},Y]\bt{w} =0$  and hence $[\tilde{X},Y] = \widetilde{[X,-\frac{i}{2}T_N]} \in \crT M$. This implies that $(M, \crT M)$ is a CR manifold. A simple computation that shows that $\crT M \subset \Ker{\ut}$ and so $\ut$ is a pseudohermitian form. Since $T \lrcorner d\ut =0$ we immediately see that $Tw=T\bt{w}=0$ and so $T=\widetilde{T_N}$. To show strict pseudoconvexity we note that $\pi_* [Y,\bt{Y}]=0$, $[Y,\bt{Y}]w=-2is$ and $[Y,\bt{Y}]\bt{w} = -2is$. Thus $[Y,\bt{Y}] = \bt{Y} -Y - iT$. Since the compact manifold $N$ is strictly pseudoconvex, this is sufficient. 

\epf

Since we have a CR structure admitting a pseudohermitian form we can easily now construct a compatible $J$ operator making $(M,J,\ut)$ a strictly pseudoconvex pseudohermitian manifold. 

\bgC{Consequences}
Viewed as a function on $M$, $w$ is CR.
\enC

The geometry of the pseudohermitian structure on $M$ is closely related to that of the the compact manifold $N$. As a starting point, we note the following interactions between the intrinsic Webster-Tanaka connection $\nabla$ for $M$ and $\nabla^N$ for $N$. 

\bgL{Simple} Suppose $X \in C^\infty(\cn{}TN)$.
\bgEn{\alph}
\etem $\nabla_Y \tilde{X} = -\frac{i}{2} \nabla_T \tilde{X}$, $\nabla_{\tilde{X}} Y =0$.
\etem $\widetilde{\nabla^N_X Z} = \nabla_{\tilde{X}} \tilde{Z}$.
\etem $\nabla_{\tilde{X}} (\pi^* \varphi) = \pi^* \nabla^N_X \varphi$.
\etem $\nabla_Y Y= -Y$ mod $(\crT M)^\top$, $\text{div } Y =-1$.
\enEn
\enL

\pf For part \rfI{a} we note that any $\tilde{X}$ is constant along any curve $\gamma$ such that $\pi_* \gamma\upp =0$. Thus $\nabla_Y \tilde{X} =- \frac{i}{2} \nabla_T \tilde{X}$. Likewise $\pi_* (Y + \frac{i}{2}T) =0$ so $Y+\frac{i}{2}T$ is constant along any curve tangent to the level sets of $w$.  Thus $\nabla_{\tilde{X}} Y = -\frac{i}{2} \nabla_{\tilde{X}} T =0$.

To see part \rfI{b} we construct an operator $\tilde{\nabla}^N: C^\infty(\cn{} TN) \times C^\infty( \cn{}TN) \to  C^\infty(  \cn{}TN ) $ by 
\[ \tilde{\nabla}^N_X Z  = \pi_* \nabla_ {\tilde{X}}  \tilde{Z}.\] 
This operator is clearly a connection on $N$. Straightforward computations then imply that it meets defining properties of the Webster-Tanaka connection and thus by uniqueness $\tilde{\nabla}^N = \nabla^N$. From \rfI{d} it follows that $\aip{\nabla_{\tilde{X}} \tilde{Z}}{Y}{} = \aip{\nabla_{\tilde{X}} \tilde{Z}}{\bt{Y}}{}=0$. Thus $\nabla_{\tilde{X}} \tilde{Z}$ annihilates $w$ and $\bt{w}$. This suffices for \rfI{b}. Part \rfI{c} is an easy consequence of \rfI{b}.

For part \rfI{d} we note that $\aip{\nabla_Y Y}{Y}{} = -\aip{Y}{\nabla_\bt{Y} Y}{}$. However $\nabla_\bt{Y} Y = [\bt{Y},Y]_{\crT} = Y$. From \rfI{a} it follows that $\text{div } Y = \ut^0 (\nabla_Y Y) =-1$.

\epf

At last, we have the background to introduce the model class of examples for the study of the $\db$-Neumann problem. Consider those domains of the form
\bgE{Model} \Omega = D \times N  \subset M\enE
where $D$ is a smoothly bounded precompact domain in \hy{2}. These domains then satisfy two important properties
\bgEn[Omega]{\Roman}
\etem $\Omega$ has no characteristic boundary points.
\etem $\Omega$ admits a smooth defining function $\rho$ depending solely on the real and imaginary parts of the CR function $w$.
\enEn
The second of these is an immediate consequence of the definitions. To see the first just note that if $Y\rho =0$ at some point $p$ then $\bt{Y}\rho =0$ at $p$, thus $d\rho =0$ at $p$. Since $\rho$ is a defining function for $\Omega$ it then follows that $Y\rho \ne 0$ on $\dO$. Thus $\dO$ is entirely non-characteristic.

\bgR{Noncompact}
In Part 2 we shall focus on the non-compact domain $D=\{|w|<1\}$ in \hy{}. In the special case when $N$ is the unit sphere the Heisenberg group \hn{2n+1} can be viewed as a completion of $M$ across the line $\{s=0\}$. For this reason we shall conduct our analysis on the set $D$ in such a way that the results apply to both the compact case and this one. 
\enR

The Kohn Laplacian depends upon the $L^2$ structure imposed by the pseudohermitian form. Recall that the volume form on $M$ is given by $dV = \ut \wedge d\ut^n$.  A simple computation thus yields that 
\bgE{Volume} dV = \frac{n}{2} dV_N \wedge \frac{i}{2} dw \wedge d\bt{w}.\enE
We shall add some flexiblity by including the presence of a weight factor $s^{\nu}$ with $\nu \geq 0$ a fixed constant, i.e.
\bgE{L2nu} \ip{\varphi}{\psi}{L^2(\Omega)} = \int\limits_\Omega \aip{\varphi}{\psi}{} s^{\nu} dV. \enE
We shall also suppose that Kohn Laplacian is defined in terms of this weighted structure, but this shall largely be suppressed from the notation.

\bgD{Components}For a $(0,q)$-form $\varphi$ on $M$ we define the transverse component of $\varphi$ by $\varphi^\bot := \bt{Y} \lrcorner \varphi$. Thus $\varphi^\bot$ is a $(0,q-1)$-form. We also introduce the tangential component of $\varphi$ by $\varphi^\top := \varphi - \ut^\bt{0} \wedge \varphi^\bot$. 
\enD

\bgL{Elementary}  Let $\varphi$ and $\psi$ be smooth $(0,q)$-forms. Then
\bgEn{\alph}
\etem $\varphi = \varphi^\top + \theta^\bt{0} \wedge \varphi^\bot$.
\etem $\aip{\varphi}{\psi}{} = \aip{\varphi^\top}{\psi^\top}{} + \aip{\varphi^\bot}{\psi^\bot}{}$.
\etem $ \left( \varphi^\bot \right)^\top = \varphi^\bot$, $\left( \varphi^\top \right)^\bot =0$.
\enEn
\enL

Only the second of these requires any work, but it follows easily from the observation that $Y$ is perpendicular to $(\crT M)^\top$ in the Levi metric.

Define an operator $\db[\top]$ by $\db[\top]:= \db - \ut^\bt{0} \wedge \nabla_\bt{Y}$. On $(0,q)$-forms it follows immediately from \rfL[PS]{Computedb} that $\db[\top] \pi^* u = \pi^* \db[N]u$. We now use this operator to split the components of the Kohn Laplacian into tangential and transverse pieces.

\bgL{Facts} Suppose $\varphi$ is a $(0,q)$-form. Then
\bgEn{\alph}
\etem $(\db[\top] \varphi)^\top = \db[\top] (\varphi^\top)$, \quad $(\db[\top] \varphi)^\bot = - \db[\top] (\varphi^\bot)$.
\etem $(\nabla_{\bt{Y}} \varphi)^\top = \nabla_{\bt{Y}} (\varphi^\top)$, \quad $(\nabla_{\bt{Y}} \varphi)^\bot = (1+\nabla_{\bt{Y}}) (\varphi^\bot)$
\etem $(\db \varphi)^\top = \db[\top] \varphi^\top$ and  $(\db \varphi)^\bot = \nabla_{\bt{Y}} (\varphi^\top) - \db[\top](\varphi^\bot).$
\etem $(\dbs \varphi)^\top = \dbs[\top] \varphi^\top + (\nabla_{\bt{Y}})^* \varphi^\bot$ and $(\dbs \varphi)^\bot =   \dbs[\top] \varphi^\bot$.
\enEn
\enL

\pf
Parts \rfI{a} and \rfI{b} are easy consequences of the observations that $\db \ut^\bt{0} = 0$ and $\nabla_\bt{Y} \ut^\bt{0} =\ut^\bt{0}$. The other parts follow easily from the definitions.

\epf

\bgL{KohnLaplacian} Define two unbounded operators on $L^2(\Lambda_\ut^{0,q}M)$ by
\[ \Pt[] \varphi := \left\{ \db[\top]\dbs[\top] + \dbs[\top]\db[\top]  + \nabla_{\bt{Y}}^* \nabla_\bt{Y} \right\}\varphi, \qquad \Pb[]\varphi := \left\{ \db[\top]\dbs[\top] + \dbs[\top]\db[\top]  + \nabla_{\bt{Y}}^* \nabla_\bt{Y} \right\} \varphi.\]
Suppose the $(0,q)$-form $\varphi \in \dom{\boxb}$. Then $\varphi^\top \in \dom{\Pt[]}$, $\varphi^\bot \in \dom{\Pb[]}$ and 
\[ \boxb \varphi = \Pt[] \varphi^\top + \ut^\bt{0} \wedge \Pb[] \varphi^\bot.\]
Furthermore both $\Pt[]$ and $\Pb[]$ preserve the tangential properties of forms, so this is an orthogonal split.
\enL

\pf
It follows easily from \rfL{Facts} that
\[
\boxb \varphi = \Pt[] \varphi^\top + [\db[\top],\nabla_\bt{Y}^*]\varphi^\bot + \ut^\bt{0} \wedge \left( \Pb[]\varphi^\bot + [\nabla_\bt{Y},\dbs[\top]]\varphi^\top \right).
\]
Now the formal adjoint of $\nabla_\bt{Y}$ is given by $\text{div }Y -\nabla_Y = -1-\nabla_Y$. By \rfL[NM]{Basic} we see $[\nabla^N_{T_N},\db[N]]=0$ from which it easily follows that $[\nabla_Y, \db[\top]]=0$. Both commutation terms thus vanish.
 
\epf

We can now use the harmonic analysis on the normal manifold $N$ to construct a partial Fourier transform for $M$. For $\nu \geq 0$ we shall denote by $\Vv$ the collection of smooth forms $\{( \sqrt{2/n}) s^{-\nu/2} \pi^* \sigma: \sigma \in \V\}$. A typical element of $\Vv$ will still be denoted by $\sigma$.  This definition is motivated by the observation that by \rfE{Volume} if $f_1,f_2 \in L^2(D)$ and $\sigma_1, \sigma_2 \in \Vv$ then
\bgE{BasicL2} \ip{f_1 \sigma_1}{f_2\sigma_2}{L^2(\Omega)} =\begin{cases} \ip{f_1}{f_2}{L^2(D)} \quad &\text{if $\sigma_1 =\sigma_2$}\\
0, \quad &\text{otherwise.}\end{cases}\enE
where the $L^2$ structure on $D$ is induced by hyperbolic metric on $\hy{2}$. The weighting by $s^{\nu}$ does not affect the adjoints of operators that act tangentially to the level sets of $w$. Thus the functions $\Gamma$, $\ul$ and $G$ are easily transferred over to the new sets $\Vv$ without change. 

We can now construct the basic decomposition that will underlie all our subsequent analysis. 

\bgD{Coefficient}
Suppose $\varphi \in L^2(\Lambda_\ut^{0,q}\Omega)$ (weighted by $\nu$). Then the tangential partial Fourier components of $\varphi$ are 
\[ \varphi^\top_\sigma (z):= \int\limits_{w=c} \aip{\varphi^\top}{\sigma}{} s^{\nu} \pi^* dV_N, \qquad \text{for $\sigma \in \Vv$.}\]  The transverse partial Fourier components are 
\[ \varphi^\bot_\sigma(z) := \int\limits_{w=c} \aip{\varphi^\bot}{\sigma}{} s^{\nu}\pi^* dV_N, \qquad \text{for $\sigma \in \Vv[q-1]$.}\]
\enD

\bgL{Decomposition}
For $\varphi \in L^2(\Lambda_\ut^{0,q}\Omega)$ each  $\varphi^\top_\sigma$, $\varphi^\bot_\sigma \in L^2(D)$ and
\[ \varphi = \sum\limits_{\Vv} \varphi^\top_\sigma \sigma  + \ut^\bt{0} \wedge \sum\limits_{\Vv[q-1]} \varphi^\bot_\sigma  \sigma.\]
Furthermore if $\psi \in L^2(\Lambda_\ut^{0,q}\Omega)$ then
\[ \ip{\varphi}{\psi}{L^2(\Omega)} = \sum\limits_{\Vv} \ip{\varphi^\top_\sigma }{\psi^\top_\sigma}{L^2(D)} + \sum\limits_{\Vv[q-1]} \ip{\varphi^\bot_\sigma}{\psi^\bot_\sigma}{L^2(D)}.\]
\enL

\pf
The proof follows easily from the observation that the original set $\V$ was an orthonormal basis for $L^2$ $(0,q)$-forms on $N$ and that each $\{w=c\}$ is isomorphic to $N$ as a pseudohermitian manifold.  

\epf

Under appropriate identifications we can write $\bt{Y} = \nW + \frac{i}{2}T$. Now $\nW s^{-\nu/2} = -\nu/2 s^{-\nu/2}$ and by unwinding the definitions we see $-i\nabla_T \sigma = \ul(\sigma) \sigma$. This implies that $\nabla_\bt{Y}  \sigma= -\left\{ (\ul(\sigma) +\nu)/2 \right\} \sigma$.

\bgD{nWs}
Define $\nWs$ and $\Ws$ as unbounded operators on $L^2(D)$ to be the maximal, closed extensions of 
\[ \nWs = \nW -\frac{\ul(\sigma)+\nu}{2}, \qquad \Ws = W +\frac{\ul(\sigma) - \nu}{2}.\]
\enD
This definition is motivated by the discussion above; for $f\in L^2(D)$, $\nabla_\bt{Y} f\sigma = (\nWs f)\sigma$, $\nabla_Y f\sigma = (\Ws f) \sigma$ where these are defined. The reader is cautioned to note that $\nWs$ and $\Ws$ are not conjugate operators.

\bgT{FourierKohn}
A $(0,q)$-form $\varphi$ is contained in $\dom{\boxb}$ if and only if all $\varphi^\top_\sigma \in \dom{ \nWs^* \nWs} $ for all $\sigma \in \Vv$ and $\varphi^\bot_\sigma \in \dom{ \nWs \nWs^*}$ for all $\sigma \in \Vv[q-1]$. Furthermore for $\varphi \in \dom{\boxb}$
\bgE{Kohn} \boxb \varphi = \sum\limits_{\Vv} \left\{ \Gm + \nWs^* \nWs \right\} \varphi^\top_\sigma \, \sigma + \ut^\bt{0} \wedge \sum\limits_{\Vv[q-1]} \left\{ \Gm + \nWs \nWs^* \right\} \varphi^\bot_\sigma \, \sigma.\enE
\enT 

\pf
For $\varphi \in \dom{\boxb}$ the expression \rfE{Kohn} for \boxb follows easily from \rfL{KohnLaplacian}. Each component of the operator defined by \rfE{Kohn} is self-adjoint as and operator on $L^2(D)$. As the operator diagonalises, it is self-adjoint as an operator on $L^2(\Omega)$. However it is a closed extension of the self-adjoint operator \boxb and so must equal \boxb exactly.   

\epf

This result can be used to improve the result of \rfL{KohnLaplacian}. We now see that $\boxb = \Pt[] + \ut^\bt{0} \wedge \Pb[]$ as unbounded operators on $L^2(\Omega)$.

We now introduce a family of Sobolev spaces for the domain $D$ which are closely related to the Folland-Stein spaces on $\Omega$.
\bgD{Sobolev}
For $\sigma \in  \Vv$ we define $\WjsD{1}$ to be the collection of all $u \in L^2(D)$ such that 
\[ \norm[2]{u}{\WjsD{1}} := \norm[2]{\G u}{L^2(D)} + \norm[2]{\Ws u}{L^2(D)} + \norm[2]{\nWs u}{L^2(D)}< \infty.\]
For $k>1$ the spaces $\WjsD{k}$ are then inductively defined as the collections of $u \in \WjsD{k-1}$ such that
\[ \norm[2]{u}{\WjsD{k}} := \norm[2]{\G u}{\WjsD{k-1}} + \norm[2]{ \Ws u}{\WjsD{k-1}} + \norm[2]{\nWs u}{\WjsD{k-1}} < \infty.\]
\enD
At each degree these spaces are pairwise equivalent  but the constants cannot be chosen to be uniform over all choices of $\sigma$. However this does yield a well-defined set $\Wo$ which is defined to be the closure of $\Cic{D}$ under any of the norms $\WjsD{1}$. This allows to define Sobolev spaces with Dirichlet boundary conditions by $\Wo[k] = \WjsD{k} \cap \Wo$.

\bgL{FourierFS}
For $(0,q)$-forms, the $k$th- order Folland-Stein norm $\norm{\cdot}{\SjO{k}}$ is equivalent to the norm defined by
\[ |\!|\!| \varphi|\!|\!|^2:= \sum\limits_{\Vv} \norm[2]{\varphi^\top_\sigma}{\WjsD{k}} + \sum\limits_{\Vv[q-1]} \norm[2]{\varphi^\top_\sigma}{\WjsD{k}}.\]
\enL
This lemma is an easy consequence of the definitions and \rfL[NM]{FollandStein} and can easily be extended to various weighted versions. In the sequel we shall use these norms interchangeably and denote both by $\norm{\cdot}{\SjO{k}}$. 

We have now reduced the study of the $\db$-Neumann problem on $\Omega$ to the analysis of two separate infinite families of operators on the domain $D$. Each family is elliptic in the interior, although only the transverse family shall prove elliptic up to the boundary. In order to obtain estimates on $\Omega$ itself we shall need uniform estimates on these families in terms of a varying family of Sobolev norms associated to the members of  \Vv. Since the family of elliptic operators under study will be closely related to these norms, this variation shall prove a boon. 
\section{Analysis on the Hyperbolic Plane}\setS{HY}

In the previous section we reduced the study of the Kohn Laplacian to understanding an infinite family of operators on the hyperbolic plane. Here we review some of the basic analysis of these types of operators. Throughout this section we suppose that $D$ is a smooth domain in \hy{2} bounded in the sense that
 \bgE{Set}
\sup\limits_{w \in D} |w| <\infty.
\enE

We shall work with the following collection of Hermitian forms. For $\ua,\ub \in \rn{}$ and sufficiently smooth functions $u$ and $v$ we define
\begin{align*}
\bQ{\ua}(v,u) &= \ip{(\ua+W)v}{(\ua+W)u}{},\\
\btQ{\ub}(v,u) &= \ip{(\ub+\nW )v}{(\ub+\nW )u}{}.
\end{align*}
Here we have employed the convention that all unlabelled norms and inner products refer to $L^2(D)$. We shall work primarily with the function space $\WjRD{1}$ which is defined as in \rfD[MC]{Sobolev} but with $\uls =\nu =0$. This space is equivalent to any individual choice of $\WjsD{1}$.

Much of our analysis in this section is based upon the following fundamental calculation. 
\bgL{Fundamental}
For any $\ua,c \in \rn{}$ the following identity holds for all $u,v \in \Wo$: 
\[\begin{split}\bQ{\ua}(v,u) &= (1-c) \bQ{(\ua-c)}(v,u) + c \btQ{(c-1-\ua)}(v,u)\\
& \qquad  + \Big( 2c\ua +c(1-c) \Big) \, \ip{v}{u}{}.\end{split}\]
\end{lemma}

\pf
Clearly $\Wo \subset \dom{\nW ^*} \cap \dom{W^*}$. Since both sides of the required identity are continuous in \WjRD{1}  it is sufficient to prove the result for $u,v\in \Cic{D}$. 

We note that $W^* =1-\nW $, so for $\ua,\ub \in \rn{}$ we compute 
\[
\bQ{\ua}(v,u) =\ip{(1+\ua-\nW )(\ua+W)v}{u}{}\]
and
\[\btQ{\ub}(v,u) = \ip{(1+{\ub}-W)(\ub+\nW )v}{u}{}.\]
Recalling that $[W,\nW ]=\nW -W$ we note
\[(1+{\ub}-W)(\ub+\nW )=  \ub(1+\ub) + \ub\nW  +(1-\ub)W -\nW W.\]
We then note that
\begin{align*}
(\ua+1-\nW )(\ua+W) & = (1-c)(\ua-c+1-\nW )(\ua-c+W) \\
& \qquad + c(c-\ua-W)(c-1-\ua+\nW ) \\
& \qquad +\ua(\ua+1) - (1-c)(\ua-c+1)(\ua-c) - c(c-\ua)(c-1-\ua)\\
 \end{align*}
and the proof easily follows from completing the square .

\epf
 
 \bgC{FundamentalEstimate}  Suppose $\ua \in \rn{}$ 
\bgEn{\alph}
\etem If $\ua \geq \frac{1}{2}$ then $2\ua \norm{u}{} \leq \bQ{\ua}(u,u)$ for $u \in \Wo$.
\etem If $-\frac{1}{2} \leq \ua \leq \frac{1}{2}$ then $
(\ua+\frac{1}{2} )^2 \norm[2]{u}{} \leq \bQ{\ua}(u,u)$ on for $u\in\Wo$.
\etem If $\ua > -\frac{1}{2}$ then for some $C>0$ the estimate $\norm[2]{u}{\WjRD{1}} \leq C \bQ{\ua}(u,u)$ holds for $u \in \Wo$.
\enEn
\end{cor}

\pf
For \rfI{a} take  $c=1$. For \rfI{b} take $c=\ua+\frac{1}{2}$. For \rfI{c} take any $0<c<1$ and apply the previous parts.

\epf

\bgC{Riesz}
Suppose $\ua > -\frac{1}{2}$ and $\phi$ is a bounded linear functional on $\Wo$. Then there exists a $u \in \Wo$ such that
\[ \bQ{\ua}(v,u) = \phi(v)\] for all $v\in \Wo$.
\enC

\pf
The Riesz representation theorem for the Hilbert space $\Wo$ (with inner product induced by $\WjRD{1}$) implies the existence of an element $z \in \Wo$ such that $\phi(v)=\ip{v}{z}{\WjRD{1}}$. By \rfCi{FundamentalEstimate}{c},  the form $\bQ{\ua}$ is strictly coercive over $\Wo$. Thus the standard Lax-Milgram Lemma can be applied to find a $u \in \Wo$ such that $\bQ{\ua}(v,u)=\ip{v}{z}{\WjRD{1}}$. 

\epf

\bgL{Shift}
For any $\ud \in \rn{}$ there is a constant $C(\ud)>0$ depending solely on $\ud$ such that
\[  \norm[2]{u}{}+ \bQ{\ua}(u,u) \leq C(\ud) \left\{ \norm[2]{u}{}+ \bQ{(\ua+\ud)}(u,u) \right\} \] for all $u \in\Wo$ and $\ua \in \rn{}$. 
\end{lemma}

\pf 
If $\ud=0$ then any $C \geq 1$ works. Suppose $\ud \ne 0$ and choose $c$ such that $1<c^2<1+\ud^{-2}$. Then
\begin{align*}
\bQ{(\ua+\ud)}(u,u)+ \norm[2]{u}{} & \geq  \bQ{\ua}(u,u) + (1+\ud^2) \norm[2]{u}{}\\
& \qquad  - 2\ud \norm{u}{} \sqrt{\bQ{\ua}(u,u)}\\
& \geq (1 +\ud^2 -c^2 \ud^2) \norm[2]{u}{} + (1 -\frac{1}{c^2})\bQ{\ua}(u,u)\\
& \geq \min\left\{ 1-\frac{1}{c^2}, 1+\ud^2(1-c^2)\right\} \left(  \bQ{\ua}(u,u) + \norm[2]{u}{} \right).
\end{align*}
This constant is strictly positive and is independent of the choice of $\ua$ and $u$.

\epf

\bgT{Precompact}
Suppose $D$ is precompact in \hy{2}.  Then for all $\ua \in \rn{}$ 
\bgEn{\alph}
\etem The bounded operator $W+\ua\colon\Wo \to L^2(D)$ is injective with closed range.
\etem The unbounded operator $\nW+\ua\colon L^2(D) \to L^2(D)$ is surjective with infinite dimensional kernel.
\enEn
\enT

\pf
When $\ua > -\frac{1}{2}$, \rfI{a} follows immediately from \rfC{FundamentalEstimate}. For $\ua \leq -\frac{1}{2}$ we first note that if $u \in \Ker{W+\ua}$ then $s^\ua u$ is holomorphic on $D$. Since $D$ is precompact and $u\in \Wo$ this implies that $s^\ua u$ is a harmonic function that vanishes on $\partial D$. Thus $u =0$ and $W+\ua$ is injective. Now apply \rfL{Shift} and \rfL{Fundamental} with $0<c<1$ to see that
\[ \norm[2]{u}{\WjRD{1}} \leq C_1  \bQ{0}(u,u) \leq C_2 \left\{ \norm[2]{u}{} + \bQ{\ua}(u,u)\right\}.\] Since $W+\ua$ is injective, the Rellich Lemma then implies that there is an estimate $\norm[2]{u}{\WjRD{1}} \leq C \bQ{\ua}(u,u)$ for $u\in \Wo$ with some constant $C>0$. For suppose there is a sequence $u_j\in \Wo$ such that $\norm{u_j}{\WjRD{1}}=1$ for all $j$ and $\bQ{\ua}(u_j,u_j) \to 0$. By the Banach-Alaoghu Theorem there is a subsequence that converges weakly in $\Wo$ to some $u \in \Wo$. By the Rellich Lemma this subsequences converges strongly to $u$ in $L^2(D)$. Thus the limit $u\in \Wo$ satisfies $\bQ{\ua}(u,u)=0$.   Injectivity of $W+\ua$ then implies that $u =0$, which contradicts the assumption that each $\norm{u_j}{\WjRD{1}}=1$. This establishes \rfI{a}.

For \rfI{b} we note that the map $v \mapsto \ip{v}{f}{}$ is bounded as a linear functional on $\Wo$. From the Rellich Lemma computation of \rfI{a} we note that the result of   \rfC{Riesz} can be extended to all $\ua \in \rn{}$ when $D$ is precompact. This implies the existence of a unique $u \in \Wo$ such that $\bQ{-1-\ua}(v,u)=\ip{f}{v}{}$ for all $v \in \Wo$. Since the formal adjoint of $W+\ua$ is  $-(\nW -1- \ua)$ we see $f= (\nW +\ua)(W-\ua-1)u$ in the sense of distributions. This proves the slightly stronger result that $(\nW+\ua)(\nW+\ua)^*$ is surjective as an unbounded map from $\Wo$ to $L^2(D)$.

The kernel of $W+\ua$ as an unbounded map on $L^2$ is given by functions of the type $s^{-\ua} \cdot h$ where $h$ is holomorphic. This set is clearly infinite dimensional.

\epf

\bgP{Domain}
For all bounded domains $D$ the $L^2$ domain of the operator $\nW^*$ is exactly $\Wo$.
\enP

\pf
Suppose $u \in \dom{\nW ^*}$ and $\ub<-\frac{1}{2}$ . Then $f = (\ub+\nW )^*u \in L^2(D)$ and  the functional $v \mapsto \ip{(W-\ub-1)v}{f}{}$ is bounded on $\Wo$. The formal adjoint of $\ub +\nW$ is given by $-(W -\ub -1)$. Since $\ub>-\frac{1}{2}$ we can apply \rfC{Riesz} to establish the existence of a unique element $z \in \Wo$ such that $\bQ{-1-\ub}(v,z)=\ip{(W-\ub-1)v}{f}{} $ for all $v \in \Wo$. This implies that \[\ip{(\ub+\nW )^* v}{ (\ub+\nW )^*(z+u)}{}=0\] for all $v \in \Wo$. We can repeat the argument of \rfi{Precompact}{b} to see that the operator $(\nW+\ub)(\nW+\ub)^*$ is surjective as an operator from $\Wo \to L^2(D)$.  Thus $z+u=0$ and so $u\in \Wo$. 

\epf

\bgT{Disc}
Suppose $D = \{|w|<1\} \subset \hy{2}$. 
\bgEn{\alph}
\etem If $\ua >-\frac{1}{2}$ then 
\bgEnn[innera]{\roman} 
\eetem The bounded operator $W+\ua\colon \Wo \to L^2(D)$ is injective and has closed range.
\eetem The unbounded operator $\nW+\ua\colon L^2(D) \to L^2(D)$ is injective with dense range, but not surjective.
\enEn
\etem If $\ua < -\frac{1}{2}$ then
\bgEnn[innerb]{\roman}
\eetem The bounded operator $W+\ua\colon \Wo \to L^2(D)$ is injective with dense range, but not surjective.
\eetem The unbounded operator $\nW+\ua\colon L^2(D) \to L^2(D)$ is surjective with infinite dimensional kernel.
\enEn
\enEn
\enT

\pf
Part \rfii{innera}{a} follows immediately from \rfC{FundamentalEstimate}. The surjectivity of part \rfii{innerb}{b} has already been established in the argument for \rfP{Domain}. The kernel consists of the products of $s^{-\ua}$ with holomorphic functions which is easily seen to be infinite dimensional.

Injectivity for \rfii{innera}{b} follows from the observation that if $u \in \Ker{\nW +\ua}$ then $s^{\ua}u$ is holomorphic. Furthermore since $u \in L^2(D)$ we can see that $s^\ua u$ extends continuously to zero across the line $\{s=0\}$. The reflection principle implies we can extend it holomorphically across the axis and then the identity theorem implies that $u$ vanishes identically.  To see that the operator is not surjective we construct the following example: Let $h_j$ be a sequence of holomorphic functions that tend uniformly to $1$ on some small disc around $i/2$ and uniformly to zero on some strip $\e<\imag{w}<2\e$. Now choose a smooth function $\phi$ that is equal to $1$ on $\imag{w}>2\e$ and $0$ on $\imag{w}<\e$ for some small $\e>0$. Then the sequence $f_j =\phi s^{-\ua} h_j$ is uniformly bounded away from $0$ in $L^2$ but $(\nW+\ua)f_j = (\nW\phi)s^{-\ua}h_j$ tends uniformly to zero.  As $\nW+\ua$ is injective this implies it does not have closed range so cannot be surjective. 

 Injectivity for \rfii{innerb}{a} follows from the arguments of \rfT{Precompact}. By \rfP{Domain} we have $W+\ua = -(\nW-1-\ua)^*$ exactly as an operator $\Wo$ to $L^2(D)$. Thus dense range follows from injectivity in \rfii{innera}{b}.

\epf

\section{Transverse Regularity}\setS{TV}

Let $D$ be a smoothly bounded domain in \hy{2}, we shall suppose that the pair $(D,\nu)$ satisfies one of the following conditions: (the results of Part 1 will focus on \rfi{C}{a}, Part 2 on \rfi{C}{b})
\bgEn[C]{\condition{C}}
\etem $D$ is precompact in \hy{2} and $\nu \geq 0$.
\etem $D= \{|w|<1\}\subset \hy{2}$ and $\nu \geq 2$.
\enEn
Throughout this section, we shall adopt the convention that unlabelled norms and inner products refer to $L^2(D)$.
 
We define a continuous function $\ddDn$ on $\hy{2}$ by setting $\ddDn[w]$ to be the hyperbolic distance from $w$ to $\dD$. We now construct  a constant $\ud(D)>0$ and smooth real valued  function $\varrho \in C^\infty(\hy{2})$  such that 
\bgEn[R]{\condition{R}}
\etem $\varrho = \ddDn$ on $\{\ddDn < 3\ud(D)\} \cap D$.
\etem $\varrho=1$ on $\{\ddDn>4\ud(D)\} \cap D$.
\etem $\varrho=0$ on $\dD$, $\varrho>0$ on $D$ and $\varrho<0$ on $\hy{2} \backslash D$.
\enEn
Since $D$ satisfies one of \rfi{C}{a} or \rfi{C}{b} this function $\varrho$ can be constructed to possess the additional important properties
\bgEnS[R]{\condition{R}}{3}
\etem $b^1 := \inf\limits_{\ddDn < 3 \ud(D)}  \min\{|W\varrho|,|\Bt{W}\varrho|\}  >0$.
\etem $B^m := \sup\limits_D  \max\limits_{j+k \leq m} |W^j \bt{W}^k \varrho| < \infty$ for all $m \geq 0$.
\enEn

If $D$ is precompact then the above conditions are satisfied by any smooth defining function. When $D=\{|w|<1\}$ then we note that $\ddDn$ is fixed under any M\"obius transformation preserving $\dD$. Thus the magnitude of $W$ is preserved by M\"obius transformations, the sizes of the various derivatives of $\ddDn$ at any point are determined purely by their values on the line $t=0$. But along $t=0$ we can explicitly compute $\ddDn(t+is)= - \ln s$. By symmetry we this line critical in the $t$-directions and so   $W\ddDn = s \partial_s (- \ln s) =-1$ along $t=0$. A simple smoothing out process then completes the construction of a suitable $\varrho$.

The extra precision is needed in the second case to allow us to construct estimates that are uniform not only across choice of $\s$ but also across a countable cover of the domain. This uniformity is essential to obtaining regularity and estimates in Part 2 of this paper. The conditions \rfi{R}{b} and \rfi{R}{c} on $\varrho$ imply that it cannot be continuously extended to the entire line $\{\imag{w} =0\}$; it would need to take the value $0$ at the boundary points at infinity while taking the value $1$ at interior points at infinity.

To avoid integrability issues when $D$ is not precompact, we shall need a refined class of smooth functions. 
\bgD{Spaces}\hfill

\bgEn{\arabic}
\etem The set $\CpD$ is the collection of all $ u\in \CiD$ such that there is some $\e>0$ with $u=0$ on $\{ \imag{w} < \e\}$.
\etem The set $\CvD$ is the collection of all $u \in \CpD$ such that $u=0$ on $\dD$.  
\enEn
\enD

For a fixed predetermined value of $q$, we shall use the notation $f(\s,u) \lesssim g(\s,u)$ on $\mathscr{W}$ to indicate the existence of a constant $C>0$ such that $f(\s,u) \leq C g(\s,u)$ for all $u \in \mathscr{W}$ and $\s \in \Vv$. Throughout this section if an explicit $\mathscr{W}$ is not mentioned we shall take it to be $\CvD$.  For convenience of reference we shall also set $ \ua=\za:=(\nu+\uls)/2-1$. The motivation behind this now being that formally
\[ \nWs^* = - \Wa\]
where $\Wa = W+\ua$.

\bgL{Prelim}
If $0 \leq q \leq n-2$ then on $\CvD$
\bgEn{\alph}
\etem $1 \lesssim \Gm + \ED{\ua}$.
\etem $\snorm{\uls}{} \lesssim \Gm + \ED{\ua}$.
\etem $\snorm{\uls}{} \lesssim \G^2$.
\etem $\G^2 \lesssim \Gm +\ED{\ua}$.
\enEn
\enL

\pf
We start by proving \rfI{a}. By \rfT[NM]{qNonZero} we see that there is a positive constant $\Go$ such that either $\Gm =0$ or $\Gm > \Go$. Thus we need only to establish a uniform estimate over those $\s$ with $\Gm=0$.  Choose $\e \in (0,\min \{\Go,\frac{1}{4}\})$ and suppose $ \uls \notin  [-\e/(n-q-1), \frac{1}{2} + \e -\nu/2]$. Then by \rfE[NM]{Bar} either $\Gm>0$ or $\ua > -\frac{1}{2} + \e$. The second of these implies $\ED{\ua} > \e^2$. Thus the only $\s$ for which we have not yet established a uniform estimate are those for which $(\Gm,\uls)\in \{0\} \times [-\e/(n-q-1), \frac{1}{2} + \e -\nu/2]$. By \rfT[NM]{qNonZero}  we see that the set of points $(\Gm,\uls)$ for $\s \in \Vv$ is discrete. Thus only finitely many can be contained in this compact set. If $\nu \geq 2$ then this set is actually empty and we are done. If $0 \leq \nu <2$ then we must be assuming   \rfi{C}{a}  and so we can use \rfC[HY]{FundamentalEstimate} to see that $\ED{\ua}>0$ for all $\s \in \Vv$. Since we have only a finite number of outstanding values of $\ua$, this is sufficient to show \rfI{a}.

Although \rfI{b} would follow from \rfI{c} and \rfI{d} we shall prove it first as it will simplify the arguments for the  later parts. We split the proof into two cases. First suppose $\ua \geq 0$ so that $\snorm{\uls}{} \leq \max \{ \nu/2-1, \uls\}$. By part \rfI{a} we need only establish an upper bound on $\uls$. Since $\ua \geq 0$  it follows easily from \rfC[HY]{FundamentalEstimate} that $\ED{\ua} \geq \frac{1}{4} (\ua+1) = \frac{\nu + \uls}{8} \geq \frac{1}{8} \uls$.  Next we suppose $\ua <0$. But then
\[ -\frac{1}{n-q-1} \Gm \leq \uls \leq 2 -\nu.\]
Part \rfI{b} then follows easily from \rfI{a}.
 
For part \rfI{c} we first note that if $q=0$ then $\uls = \Bt{\Gamma}(\s) - \Gm$ and so $\snorm{\uls}{} \leq \Bt{\Gamma}(\s) + \Gm \leq \G^2$.  If $1 \leq q \leq n-2$ the result follows immediately from the observation that $1 + \boxb[N]$ is an isomorphism from $S^2(N)$ to $L^2(N)$ while $\T[N]$ is a bounded operator between the same spaces.

For part \rfI{d} we note that when $q>0$ we easily have uniformity over those $\s$ such that $\Gm \geq \Go$. The result then follows from \rfI{a}. When $q=0$ we need to work a little harder. However we note from \rfE[NM]{Bar} that $\G = 1 + (n-1)\uls + 2\Gm$ and each piece can be uniformly bounded using previous parts.

\epf

\bgL[Basic Estimate]{BasicEstimate}
Suppose $0 \leq q \leq n-2$ and $(D,\nu)$ satisfies either \rfi{C}{a} or \rfi{C}{b}. Then on \Wo \[ \norm[2]{u}{\WjsD{1}} \lesssim  \bQ{\ua} \big(u,u\big) + \Gm \norm[2]{u}{}.\] 
\enL

\pf
By density it is sufficient to show uniformity on \CvD. Recall that by definition $\ED{\ua} \norm[2]{u}{} \lesssim \bQ{\ua}(u,u)$ on \Wo. Then by \rfL{Prelim} \rfi{Prelim}{c} we see $\G^2 \lesssim \Gm + \ED{\ua}$ on \CvD. Now for $u \in \CvD$ we can apply \rfL[HY]{Shift} to see
\begin{align*}
\bQ{(\uls-\nu)/2}(u,u) + \norm[2]{u}{} &\leq C(\ua - (\uls+\nu)/2 ) \left\{ \bQ{\ua}(u,u) +  \norm[2]{u}{} \right\}\\& \leq C(\nu -1) \left\{ \bQ{\ua}(u,u) +  \norm[2]{u}{}  \right\}\\
& \lesssim \Gm  \norm[2]{u}{}  + \bQ{\ua}(u,u)
\end{align*}
where the last line follows from \rfL{Prelim} \rfi{Prelim}{a}. Using \rfL[HY]{Fundamental} we see that
\begin{align*}
\btQ{-(\nu+\uls)/2}(u,u)& = \bQ{\ua+1}(u,u) - 2(\ua+1) \norm[2]{u}{} \\
&\leq \bQ{\ua+1}(u,u) +( 1- \uls) \norm[2]{u}{} \\
& \leq C(-1) \left\{ \bQ{\ua}(u,u) +  \norm[2]{u}{}  \right\}  + \snorm{\uls}{}  \norm[2]{u}{} \\
& \lesssim \Gm  \norm[2]{u}{}  + \bQ{\ua}(u,u).
\end{align*}
Combining these results establishes the lemma.

\epf

\bgL{Injective}
If $0 \leq q \leq n-2$ and the pair $(D,\nu)$ satisfies either \rfi{C}{a} or \rfi{C}{b} then 
\bgEn{\alph}
\etem $\Pb$ is injective.
\etem $\dom{\Pb} \subset \Wo$.
\etem If $u \in \CvD$ then $\Gm  \ip{v}{u}{} + \bQ{\ua}(v,u) = \ip{v}{\Pb u}{}$ for all $v \in \CvD$.
\etem There exists $C>0$ such that 
\[ \norm{u}{\WjsD{1}} \leq   C \norm{\Pb u}{}\] for all $u \in \dom{\Pb}$ and $\s \in \Vv$. 
\etem $\Pb$ has closed range.
\enEn
\enL
\pf
Fix $\s \in \Vv$. The formal adjoint of $\nWs$ is $1-\nu/2-\Ws$. Define a second order differential operator by $R=\Gm-\nWs \Wa$. 

From the Basic Estimate it follows that for all $f \in L^2(D)$ there is some $u \in \Wo$ such that \[\Gm \ip{v}{u}{} + \bQ{\ua}(v,u) = \ip{v}{f}{}\]
for all $v \in \Wo$. This implies $R u =f$ in the sense of distributions.

If $\Gm \geq \Go$ then clearly $\Pb$ is injective. Suppose $\Gm =0$. Then by \rfE[NM]{Bar} $\uls \geq 0$. Thus by \rfT[HY]{Precompact} if \rfi{C}{a} holds, or \rfT[HY]{Disc} if \rfi{C}{b} does, we see that $\nWs$ is surjective as an unbounded operator on $L^2(D)$. If $u \in \Ker{\Pb}$ then
 \[0=\ip{u}{\Pb u}{}=\norm[2]{\nWs^*u}{}.\] Thus $u \in \Ker{\nWs^*}$. But $\nWs$ is surjective. Therefore $u=0$ and we have established rfI{a}. 

Fix $u \in \dom{\Pb}$ and set $f=\Pb u \in L^2(D)$. Solve the equation $Rv=f$ for $v \in \Wo$. Since the formal adjoint of $\nWs$ is $-\Wa$ we see $\Wo \cap \dom{R} \subset \dom{\Pb}$. Therefore $\Pb v=R v=f=\Pb u$. But $\Pb$ is injective, so $u=v \in \Wo$ and \rfI{a} holds.

Part \rfI{c} follows easily from injectivity. Part  \rfI{d} is now an immediate consequence of the Basic Estimate and implies \rfI{e}. 

\epf

\bgD{Control}
A family of differential operators $\{ \mE \} _{\s \in \V}$ is $\s$-uniform of order $k$ if each $\mE$ can be written as a polynomial of degree $\leq k$ in $\G$, $\Ws$ and $\nWs$ with coefficients that are smooth functions independent of $\s$. 
\enD
It follows immediately that for a $\s$-uniform family of order $k$ there is a constant $C>0$ such that
\[ \norm{\mE u}{\Wjs{m}} \leq C \norm{ u}{\Wjs{m+k}}\] for all $u \in \Wjs{m}$, $m \geq 0$ and $\s \in \Vv$. We shall use the notation $\Lk{k}$ to indicate an element a generic $\s$-uniform family when the precise is not important. Unwinding these definitions, we obtain the important commutation relation
\bgE{Commutation}
[ \Lk{1} ,\Lk{1}] = \Lk{1} + \uls \Lk{0}.
\enE
It should be stressed here that the family $\{\uls\}$ is $\s$-uniform of order $2$. Furthermore it is easy to see that the formal adjoint of any $\Lk{1}$ is again  $\s$-uniform of order $1$.

We shall frequently need to make the following extra suppositions about a $\s$-uniform family $\{\mE\}$ and its interaction with a pair of cut-off functions $\xi$, $\zeta$:
\bgEn[e]{\condition{e}}
\etem $\{ \mE \}$ is $\s$-uniform of order $1$.
\etem $\xi \subset \zeta$, i.e. $\zeta =1$ on the support of $\xi$.
\etem $ \mE \xi$ maps \CvD to \CvD
\etem $ \CpD \subset \dom{ [\mE \xi]^*}$
\enEn

\bgL{GEstimate}
Suppose $u \in \CvD$, and $\xi$ and $\zeta$ are cut-off functions satisfying \rfi{e}{b}. Then
\[
 \norm{ \G \xi u}{\WjsD{1}} \lesssim \norm{\zeta \Pb u}{} + \norm{ \zeta u}{\WjsD{1}}.
\]
\enL

\pf We start by computing that
\bgE{GEstimate}
\begin{split}
 \Gm \norm[2]{ \G \xi u}{} + &\bQ{\ua}( \G \xi u, \G \xi u) = \Gm \ip{ \xi \G^2 \xi u}{u}{} + \bQ{\ua}( \xi \G^2 \xi u,u)\\
 & \quad+ \ip{\Wa\G \xi u}{ [W, \xi] \G \zeta u}{} + \ip{ [\xi,W] \G \xi u}{ \Wa \zeta u}{}\\
 & = \ip{\xi \G^2 \xi u}{\zeta f}{} + \ip{\Lk{1} \G \xi u}{ \Lk{0} \G \zeta u}{} + \ip{\Lk{0} \G \xi u}{ \Lk{1}\zeta u}{}\\
&=  \ip{\Lk{1} \G \xi u}{ \zeta f}{} + \ip{\Lk{1} \G \xi u}{ \Lk{0} \G \zeta u}{} + \ip{\Lk{0} \G \xi u}{\Lk{1} \zeta u}{}
\end{split}
\enE
Then the Basic Estimate implies
\[ \norm[2]{\G \xi u}{\WjsD{1}} \lesssim \norm{\G \xi u}{\WjsD{1}} \left\{ \norm{\zeta \Pb u}{} +\norm{\zeta u}{\WjsD{1}} \right\}. \]
The result follows easily. 

\epf

\bgL{KeyStep}
Suppose the triple $\mE, \xi, \zeta$ satisfy \rfi{e}{a} -\rfi{e}{d}. Then
\[ \norm{\mE \xi u}{\WjsD{1}} \lesssim \norm{\zeta \Pb u}{} + \norm{\zeta u}{\WjsD{1}} \] on $\CvD$. Furthermore the uniform bounding constant can be chosen to depend solely on $\zeta$, $\xi$ and their derivatives.
\enL

\pf
Suppose $v \in \CvD$ and $\xi \subset \zeta\upp \subset \zeta$. We then compute
\begin{align*}
\Gm \ip{u}{\mE \xi u}{} +\bQ{\ua}\big( v, \mE \xi u\big) &= \Gm \ip{ \xi \mE^* v}{u}{} + \ip{ \Wa v}{\mE \xi  \Wa u}{}\\
& \qquad  + \ip{\Wa v}{[\mE \xi,W]u}{}\\
&= \Gm \ip{\xi \mE^* v}{u}{}  + \bQ{\ua}\big(\xi \mE^* v, u\big) \\
&\qquad + \ip{ [\xi \mE^*,W]v}{ \Wa u}{} + \ip{\Wa v}{[\mE \xi,W]u}{}\\
&= \ip{\xi \mE^* v}{\zeta \Pb u}{} + \ip{ (\Lk{1}+\uls\Lk{0}) \xi v}{\Lk{1} \zeta u}{}\\
& \qquad + \ip{\Lk{1} v}{ (\Lk{1} + \uls\Lk{0})\zeta\upp u}{}.
\end{align*}
Since $[\mE,W]$ and $[\mE^*,W]$ are both $\s$-uniform of order $1$ it is further clear that the bounding constants associated to the $\Lk{1}$ and $\Lk{0}$ terms can be controlled by a constant multiple of $\sup \snorm{\xi}{\WjR{1}}$.

Now set $v = \mE \xi u$. Using the Basic Estimate combined with the observation that  from  \rfL{Prelim} $1 + \ulsq \lesssim \G$, we obtain
\[ \norm[2]{\mE \xi u}{\WjsD{1}} \lesssim \norm{\mE \xi u}{\WjsD{1}}  \left\{ \norm{\zeta \Pb u}{}  + \norm{\G \zeta\upp u}{\WjsD{1}} \right\}.\] 
The result now follows from \rfL{GEstimate}. The restriction on the bounding constant follows from careful examination of the commutation terms.

\epf

\bgL{Estimates}
Suppose the pair $(D,\nu)$ satisfies \rfi{C}{a} or \rfi{C}{b}. Let $\xi$ and $\zeta$ be cut-off functions such that $\xi \subset \zeta$ and one of the following holds
\bgEn{\condition{$\xi$}}
\etem $\zeta$ is  supported in the interior of $D$.
\etem $\zeta$ is supported in the region $\{\ddDn < \ud(D)\}$.
\enEn
Then
\[ \norm{\xi u}{\WjsD{k+2}} \lesssim \norm{\zeta \Pb u}{\WjsD{k}} + \norm{\zeta u}{\WjsD{1}}.\]
Furthermore the bounding constants can be chosen to depend solely on the magnitudes of $\xi$, $\zeta$ and their derivatives.
\enL
\resetclaims
\pf

The proof will run as a delicate induction argument. The first step is to establish the result when $k=0$. If $\zeta$ satisfies \rfI{a} then it is easy to check that if $\mE$ equals either $\Ws$ or $\nWs$ the conditions of \rfi{e}{a}-\rfi{e}{d} are satisfied. The case $k=0$ follows immediately from  \rfL{GEstimate} and \rfL{KeyStep}. 

If instead $\zeta$ satisfies \rfI{b} the argument is a little trickier as we must be more restrictive in order to satisfy the conditions \rfi{e}{a}-\rfi{e}{d}. We set
\bgE{DefE}
\mE =( \nW \varrho) \Ws - (W\varrho)\nWs.
\enE
We note that $\mE \varrho = \real{W\varrho} \uls\varrho + i \imag{W\varrho}\nu \varrho$ which vanishes on $\dD$. Thus $\mE$ preserves $\CvD$ and $\CpD$ is in the domain of its adjoint. Elementary commutation computations reveal that $\mE$, $\xi$  and $\zeta$ satisfy \rfi{e}{a}-\rfi{e}{d}. Thus we can apply \rfL{KeyStep} to see that
\[ \norm{\mE \xi u}{\WjsD{1}} \lesssim \norm{\zeta \Pb u}{} + \norm{\zeta u}{\WjsD{1}}.\]
We can explicitly compute
\[ [\nWs,\mE] = (\nW \nW \varrho) \Ws - (\nW W \varrho) \nWs  + (\nW \varrho) (\Ws -\nWs -\uls).\]
Therefore \bgE{Switch}
\begin{split}
(\nW \varrho) \Ws^2 &=  \left\{ (\mE + (W\varrho)\nWs \right\} \Ws \\
&= \mE \Ws -(W \varrho)(\Pb -\Gm +(\nu-1) \nWs\\
& = \Ws \mE - (W \varrho)\Pb + \Lk{1} + \G^2 \Lk{0}
\end{split}
\enE
Since we are assuming \rfi{Estimates}{b} we see from \rfi{R}{d} that $\snorm{\nW \varrho}{}$ is bounded below on $\text{supp}(\zeta)$. This implies
\[ \norm{\Ws^2 \xi u}{} \lesssim \norm{\zeta \Pb u}{} + \norm{\zeta u}{\WjsD{1}}.\] If $D$ is not precompact then the uniform bounds on the derivatives of $\varrho$ are needed for this.
A very similar argument gives control over the $\nWs^2$ derivatives. All others can be deduced from \rfL{GEstimate}. Thus we have established the result when $k=0$. 

We now proceed by induction supposing we established the result at all stages prior to $k$. Note that since $\G$ acts by constant multiplication, we easily obtain
\bgE{GBit}
\begin{split}
\norm{\G \xi u}{\WjsD{k+1}} &\lesssim \norm{\zeta \Pb u}{\Wjs{D}{k}} + \norm{\zeta\upp \G u}{\WjsD{1}}\\
& \lesssim  \norm{\zeta \Pb u}{\WjsD{k}}  + \norm{\zeta u}{\WjsD{1}}
\end{split}
\enE
where $\xi \subset \zeta\upp \subset \zeta$.

 Now suppose $\zeta$ satisfies \rfi{Estimates}{b} and choose $\xi\upp$, $\zeta\upp$ and $\zeta\dupp$  so that $\xi \subset \xi\upp \subset \zeta\upp \subset \zeta\dupp \subset \zeta$. Apply \rfL{KeyStep} with $\mE$ as in \rfE{DefE} to establish
\begin{align*}
\norm{\mE^{k+1} \xi  u}{\WjsD{1}} &\lesssim \norm{\mE \xi\upp \mE^{k} \xi u}{\WjsD{1}}\\
 &\lesssim \norm{ \zeta\upp \Pb \mE^{k} \xi u}{} + \norm{\zeta\upp \mE^{k} \xi u}{\WjsD{1}}\\
&= \norm{\zeta\upp \mE^{k} \xi \Pb u}{} + \norm{\zeta\upp [\Pb,\mE^k \xi  ]u}{} + \norm{\zeta\upp \mE^{k} \xi  u}{\WjsD{1}}\\
&= \norm{\Lk{k} \zeta\dupp \Pb u}{} + \norm{ (\Lk{k+1} + \uls\Lk{k} ) \zeta\dupp u}{} + \norm{ \zeta\dupp u}{\WjsD{k+1}}\\
& \lesssim \norm{\zeta \Pb u}{\WjsD{k}} + \norm{ \zeta u}{\WjsD{1}}
\end{align*}
where the last line follows from \rfE{GBit} and induction. We can now repeatedly apply \rfE{Switch} together with commutation arguments together with similar results about $\nWs^2$ to establish the general result. Again if $D$ is not precompact, the uniform bounds on the derivatives of $\varrho$ are needed at this stage.

The argument when we assume \rfi{Estimates}{a} instead of \rfi{Estimates}{b} is similar but simpler as instead of needing to use $\mE^k$ we instead use operators of type \Lk{k}. 

In both  cases, the restriction of the bounding constants follows from careful examination of the commutation terms at each stage.

\epf

\bgL{cut-off}
Suppose that $\psi$ is an even function in $\Cic{\rn{}}$. Set $\xi_p = \psi (\dist[]{\cdot}{p}^2)$ where $\dist[]{\cdot}{p}$ is the hyperbolic distance from $d$. Then each $\xi_p$ is smooth and for each $m \geq 0$
\[ \sup \max\limits_{j+k=m} \snorm{W^j \nW^k \xi_p}{}\] is finite and independent of $p$.
\enL

The lemma follows from routine computations involving the transitive action of the group of M\"obius transformations on \hy{2}. The key observation is that for any M\"obius transformation $\phi$ there is some $a \in \cn{}$ with $\snorm{a}{}=1$ such that $\phi_* W =a W$.

\bgT{Local}
Suppose $0 \leq q \leq n-2$ and $(D,\nu)$ satisfies one of \rfi{C}{a} or \rfi{C}{b}. Then if $f \in L^2(D)$ there exists a unique $u \in \dom{\Pb}$ such that $\Pb u=f$. Furthermore if $\zeta f \in \WjsD{k}$ for some cut-off function $\zeta$ satisfying one of \rfi{Estimates}{a} or \rfi{Estimates}{b} then for any $\xi \subset \zeta$ we have  $\xi u \in \WjsD{k+2}$ and a uniform estimate
\[ \norm{\xi u}{\WjsD{k+2}} \lesssim \norm{\zeta f}{\WjsD{k}} + \norm{\zeta u}{}.\]
\enT

\pf
From the Basic Estimate we see that the self-adjoint operator \Pb is injective and has closed range. Standard analysis of unbounded operators then implies that \Pb is surjective. Thus we can solve $\Pb u =f$ for $u \in \dom{\Pb} \subset \Wo$. The regularity gain follows from the standard local regularity theory of elliptic operators.

To establish the estimates we note that $\CvD$ is dense in $ \Wo[k+2]$ (\rfL[DY]{Vanish}) . Choose $\zeta \upp$ such that $\xi \subset \zeta\upp \subset \zeta$ and construct a sequence $u_j$ in $\CvD$ such that $u_j \to \zeta u$ in $\WjsD{k+2}$. Now $\Pb$ is continuous from $ \Wo[k+2]$ to $\WjsD{k}$ so $\zeta\upp \Pb u_j \to \zeta\upp \Pb \zeta u = \zeta\upp f$. By \rfL{Estimates} there are estimates (uniform over $\s$ and $j$) of the form
\[ \norm{\xi u_j}{\WjsD{k+2}} \lesssim \norm{\zeta\upp \Pb u_j}{\WjsD{k}} + \norm{\zeta\upp u_j}{\WjsD{1}}.\]
Letting $j \to \infty$ then establishes the estimate
\bgE{Almost}
\norm{\xi u}{\WjsD{k+2}} \lesssim \norm{\zeta f}{\WjsD{k}} + \norm{\zeta u}{\WjsD{1}}.\enE
To obtain the stronger estimate, we run through the argument of \rfL{GEstimate} using $\xi u$ instead of $\xi \G u$. After an integration by parts on the final term of \rfE{GEstimate}, we obtain the estimate
\[ \norm{\xi u}{\WjsD{1}} \lesssim \norm{\xi f}{} + \norm{\zeta u}{}.\]
We can replace  $\zeta$ by some $\zeta\upp \subset \zeta$ in  \rfE{Almost} and apply this estimate to complete the proof. 

\epf

\bgT{Transverse}
Suppose $0 \leq q \leq n-2$ and $(D,\nu)$ satisfies one of \rfi{C}{a} or \rfi{C}{b}. Then for each $\s \in \Vv$ and $k\geq 0$ the operator  
\[ \Pb: \Wo[k+2] \to \WjsD{k}\]  is an isomorphism.  Furthermore the bounding constants are uniform over all choices of $\s \in \Vv$.
\enT 

\pf
Suppose $f \in \WjsD{k}$, then by \rfT{Local} there is a unique $u \in \dom{\Pb}$ such that $\Pb =f$. Since $\Pb$ is continuous between these spaces we just need to establish continuity of its inverse.

Zorn's Lemma implies the existence of  a maximal cover collection of points $\{p_j\} \subset D$ such that  the open hyperbolic balls of radius $\e/2$ about each $p_j$ are all disjoint. Clearly the balls $B_\e(p_j)$ form an open cover of $D$. Choose $\ud > \e$ and fix a $p_j$. If $B_\ud(p_j) \cap B_\ud(p_k) \ne \emptyset$ then $B_{\e/2}(p_k) \subset B_{(2\ud+\e)}(p_j)$. This last ball has finite volume so can contain at most finitely many disjoint balls of radius $\e/2$. With the number bounded by the ratio of the volumes of the balls of radius $2\ud+\e$ and $\e/2$. This ratio is independent of where the balls are centred. Thus the collection $\{B_\ud(p_j)\}$ is uniformly locally finite.

Choose $0<\e<\ud$ to be sufficiently small that any cut-off function supported in some $B_\ud(p_j)$ satisfies either \rfi{Estimates}{a} or \rfi{Estimates}{b}. For each $j$ choose cut-off functions $\xi_j \subset \zeta_j$ that depend solely on the hyperbolic distance to $p_j$ such that $\zeta_j$ is supported in $B_\ud(p_j)$ and $\xi_j =1$ on $B_\e(p_j)$. By \rfL{cut-off} these can be chosen so that the sup norms of the magnitude of the cut-off functions and their derivatives are independent of $j$. The bounding estimates of \rfL{Estimates}  when applied to this cover can then be chosen to be independent of $j$. Thus we see 
\begin{align*}
\norm{u}{\WjsD{k+2}} &\lesssim \sum\limits_j \norm{\xi_j u}{\WjsD{k+2}} \\
& \lesssim \sum_j \norm{\zeta_j \Pb u}{\WjsD{k}} + \norm{\zeta_j u}{\WjsD{1}}\\
& \lesssim \norm{f}{\WjsD{k}} + \norm{u}{\WjsD{1}}
\end{align*}
where the last line follows from uniform local finiteness. We then use the Basic Estimate to bound $\norm{u}{\WjsD{1}}$ uniformly by $\norm{f}{}$ and the proof is complete.

\epf

\bgC{Domain}
Under the conditions of \rfT{Transverse}, $\dom{\Pb} = \Wo[2]$.
\enC

\section{Tangential Estimates}\setS{TG}

The tangential components of the partial Fourier transformation are more problematic. The domain of each $\Pt$  is difficult to describe concretely and certainly contains elements not in $\Wo$. The operators $\Pt$ do not possess  strictly coercive Dirichlet forms. Indeed, the following example shows that we cannot expect a full gain of local regularity at the boundary for solutions to $\Pt u=f$ even when $D$ is precompact.

\bgX{NegRegularity}
Let $D$ be the domain $\{\snorm{w-2i}{} < 1\}$ and set $z_0 \in \partial D$ to be the point $z_0=1+2i$. Fix $q=1$ and $\nu=0$. Fix any $\s \in \V[1]$ such that $\Gm >0$. Such a $\s$ must exist by subellipticity of $\boxb[M]$.  Now let $u$ be the function $u =s^{\uls/2} \log (w-z_0)$ with the logarithm defined by cutting to the right of $z_0$. Then $u \in L^2(D)$ and $\nWs u =0$. Therefore $\Pt u = \Gm u \in L^2(D)$. However $\Ws u =s^{\uls/2}   \frac{2is}{w-z_0} + \uls u$ has a singularity at $z_0$ which is not $L^2$-integrable.
\enX
This example shows that we cannot expect a full gain of two Folland-Stein derivatives up to the boundary on the tangential component. However it does indicate some of the possibilities we might hope for. For instance, we may expect regularity for the derivative $\nWs u$ or some degenerately weighted regularity of $u$. 
 
\bgL{Basic}
Suppose $1\leq q \leq n-2$ and the pair $(D,\nu)$ satisfies either \rfi[TV]{C}{a} or \rfi[TV]{C}{b}. Then for any $\s \in \Vv$ the operator $\Pt$ has closed range as unbounded operator on $L^2(D)$. Furthermore
\bgEn{\alph}
\etem If $\Gm >0$ then $\Pt$ is a bijection from $\dom{\Pt}$ to $L^2(D)$ and for all $u \in \dom{\Pt}$
\[ \norm{u}{} \leq (\Go)^{-1} \norm{\Pt}{}.\]
\etem If $\Gm=0$ and $f  \in \left\{ \Ker{\Pt} \right\}^\bot$ there exists a unique $u \in \dom{\Pt} \cap \WjsD{2} \cap \left\{ \Ker{\Pt} \right\}^\bot$ such that $\Pt u =f$. In addition if $f \in \WjsD{k}$ then $u \in \WjsD{k+2}$ and there is a constant $C$ independent of $\s$ and $f$ such that 
\[ \norm{u}{\WjsD{k+2}} \leq C \norm{f}{\WjsD{k}}.\]
\enEn
\enL

\pf 
First we assume $\Gm>0$, thus $\Gm \geq \Go$. Then we immediately see
\[ \norm{u}{} \leq (\Go)^{-1} \norm{\Pt u}{}\] for all $u \in \dom{\Pt}$. Standard functional analysis arguments then imply that $\Pt$ is injective and has closed range. Since $\Pt$ is self-adjoint this also yields that $\Pt$ is surjective.

Now suppose $\Gm =0$.   Then $\Ker{\Pt} = \Ker{\nWs}$. Choose $f \in \left\{ \Ker{\Pt} \right\}^\bot$. Recall that formally $\nWs^* =- \Wa$. From the work of \rfS{HY} we can find $g \in \Wo$ such that $\bQ{\ua}(v,g)=-\ip{\Wa v}{f}{}$ for all $v \in \Wo$. The Basic Estimate of \rfS{TV} implies that $\norm{g}{\WjsD{1}} \lesssim \norm{f}{}$.  Since $\Gm =0$ we can apply \rfT[TV]{Transverse} to find $z \in \Wo[3]$ such that $\nWs \nWs^* z = g$. Set $u = \nWs^* z$. Then $u \bot \Ker{\nWs}$ and $\bQ{\ua}(v,\nWs u) = -\ip{\Wa v}{f}{}$ for all $v \in \Wo$. Since $\dom{\nWs^*}=\Wo$,  this is sufficient to show that $\Pt  u -f $ is orthogonal $\rng{\nWs^*}$. But both $\Pt u$ and $f$ are elements of $\Ker{\nWs}^\bot \supseteq \rng{\nWs^*}$. Thus $\Pt u = f$. The regularity theory and estimates of the previous section  now yield 
\[\norm{u}{\WjsD{k+2}} \lesssim \norm{z}{\WjsD{k+3}} \lesssim \norm{g}{\WjsD{k+1}} \lesssim \norm{f}{\WjsD{k}}.\]
Here it should be noted that the uniformity is only over those $\s$ such that $\Gm =0$. That $\Pt$ has closed range follows immediately from this estimate with $k=0$.

\epf

In the case $\Gm =0$ we can go further and immediately produce local estimates. From a simple application of \rfT[TV]{Local} we see
\bgE{Locala}
\norm{\xi u}{\WjsD{k+2}} \lesssim \norm{\xi\upp z}{\WjsD{k+3}} \lesssim \norm{\zeta\upp g}{\WjsD{k+1}} + \norm{\zeta\upp z}{} \lesssim \norm{\zeta g}{\WjsD{k+1}} + \norm{\zeta u}{}.
\enE
Choose cut-off functions $\xi \subset \xi\upp \subset \zeta\upp \subset \zeta$. Running through arguments almost identical in nature to those of \rfS{TV} we  get
\[ \norm{\xi g}{\WjsD{k+1}} \lesssim \norm{\zeta f}{\WjsD{k}} + \norm{\zeta g}{}.\]
The gain of only a single derivative is due to the presence of a derivative in the right hand side of the equation defining $g$. Combining, we get the following local estimates when $\Gm=0$,
\bgE{Local}
\norm{\xi u}{\WjsD{k+2}} \lesssim \norm{\zeta \Pt u}{} + \norm{\zeta u}{} + \norm{\zeta \nWs u}{}.
\enE

\bgC{Cohomology}
Suppose $D$ is precompact. If $\Gm \geq \Go$ then $\Pt $ is injective, otherwise $\Ker{\Pt}$ is infinite dimensional. 
\enC

\pf
The first part is a trivial consequence of \rfL{Basic}. For the second we note that if $\Gm =0$ then $u \in \Ker{\Pt}$ if and only if $s^{-\uls/2}u$ is holomorphic. Since $D$ is precompact, the space of $L^2$ functions that are holomorphic on $D$ is infinite dimensional.

\epf

\bgL{WEstimates}
Suppose $1 \leq q \leq n-2$ and $(D,\nu)$ satisfies \rfi[TV]{C}{a} or \rfi[TV]{C}{b}. If $\Pt u =f \in \WjsD{k}$ for $k \geq 0$ then $\nWs u \in \Wo[k+1] $. Furthermore there is a uniform estimate
\[ \norm{\nWs u}{\WjsD{k+1}} \lesssim \norm{f}{\WjsD{k}}.\]
\enL

\pf
If $k>0$ then we note that $\nWs f= \nWs \Pt u = \Pb \nWs u$. As $\nWs f \in \WjsD{k-1}$  we can apply \rfT[TV]{Transverse} to establish the result.

If $k=0$ we instead note that for all $v \in \Wo$ we have $\Gm \ip{v}{\nWs u}{} +\bQ{\ua}(v,\nWs u)=\ip{-\Wa v}{f}{}$. Now $\nWs u \in \dom{\nWs^*} = \Wo$ so we can apply the Basic Estimate of \rfS{TV} to see that 
\[ \norm[2]{\nWs u}{\WjsD{1}} \lesssim \norm{\nWs u}{\WjsD{1}} \norm{f}{}.\]  

\epf

\bgC{Exact}
Under the same conditions as \rfL{WEstimates}, if $u \in \left\{ \Ker{\Pt} \right\}^\bot $ and  $\Pt u \in \WjsD{k}$ then $u \in \WjsD{k}$ with a uniform estimate
\[ \norm{\G^2 u}{\WjsD{k}} \lesssim \norm{\Pt u}{\WjsD{k}}.\]
\enC

\pf
If $\Gm =0$ then in \rfL{Basic} we have already established a vastly improved result. Namely $u \in \WjsD{k+2}$ with estimate. If $\Gm \geq \Go$ then by the previous lemma \begin{align*}
 \Gm u &= \Pt u - \nWs^* \nWs u  = \Pt u +\Wa \nWs\\
 & = \Pt u + \Ws \nWs u + (\nu-1)\nWs u \in \WjsD{k}.
 \end{align*} 
Furthermore since $\nu$ is a fixed constant and $\Ws$ is a $\s$-uniform operators of order $1$, we also obtain a uniform estimate for $\Gm u$. Since we currently assuming $q>0$ and $\Gm>0$ we have the estimate $\G^2 \leq \left(1+(\Go)^{-1} \right) \Gm$. This is sufficient to complete the proof.

\epf

Combined with the local estimates of \rfE{Local}, this quickly yields the local estimates
\bgC{Local}
For cut-off functions $\xi$, $\zeta$ satisfying \rfi[TV]{Estimates}{a} and \rfi[TV]{Estimates}{b} the following local estimate holds uniformly for $u \in \{\Ker{\Pt}\}^\bot$
\[\norm{\xi \G^2 u}{\WjsD{k}} \lesssim \norm{\zeta \Pt u}{\WjsD{k}} + \norm{\zeta u}{} + \norm{\zeta \nWs u}{}.\]
\enC

\bgR{reg}
The results we have already shown are already sufficient to establish a sharp regularity theorem for each $\Pt$ with uniform estimates. Namely, we have that for $k \geq 0$
\[ \Pt: \{ u \in \WjsD{k}: \nWs u \in \Wo[k+1] \} \cap \left\{ \Ker{\Pt} \right\}^\bot \longrightarrow \WjsD{k} \cap \left\{ \Ker{\Pt} \right\}^\bot\]
is an isomorphism. 

This initially looks a little weaker than the result of \rfL{Basic} \rfi{Basic}{b} when $\Gm =0$. But it is easy to check that $\Pt$ is continuous between these spaces. The existence and uniqueness of solutions then implies 
\[\{ u \in \WjsD{k}: \nWs u \in \Wo[k+1] \} \cap \left\{ \Ker{\Pt} \right\}^\bot = \{ u \in \WjsD{k+2}: \nWs u \in \Wo[k+1]\} \cap \left\{ \Ker{\Pt} \right\}^\bot \]
when $\Gm=0$.
\enR

\section{Weighted Tangential Estimates}\setS{WR}

The results of the previous section are not that satisfying; they do not answer the question of whether we obtain any form of gain of regularity generic directions.  By \rfX[TG]{NegRegularity}, we know that we cannot expect any gain of Folland-Stein regularity for the derivatives $\Ws u$, however we can establish weighted estimates. 

Throughout this section we shall again adopt the convention that unlabelled norms and inner products are taken to be with respect to $L^2(D)$.


\bgD{Weighted}
The spaces $\R{k} \Wjs[(U)]{j}$ are defined inductively by
\[ \R{0} \Wjs[(U)]{j} = \Wjs[(U)]{j}\] with equivalent norms and
\[ \R{k+1} \Wjs[(U)]{j} = \left\{ u \in \R{k}\Wjs[(U)]{j}: \varrho \Ws, \varrho \nWs \in \R{k}\Wjs[(U)]{j} \right\}.\]
The norms for $k>0$ are inductively defined by
\[ \norm[2]{u}{ \R{k+1}\Wjs[U]{j} } = \norm[2]{ \varrho \G u}{\R{k}\Wjs[(U)]{j}} +  \norm[2]{ \varrho \Ws u}{\R{k}\Wjs[(U)]{j}}+ \norm[2]{ \nWs u}{\R{k}\Wjs[(U)]{j}}.\]
\enD

\bgD{Local}
The space $\WlD{k}$ is defined to be the collection of all $u \in L^2(D)$ such that for all precompact sets $V \subset \hy{2}$ we have $ u \in \Wjs[(V \cap D)]{k}$. We shall make analogous definitions for other function spaces.
\enD
It is worth noting that this definition implies local differentiability at the boundary as well as in the interior of $D$. For precompact domains $D$ we easily see $\WlD{k}=\WjsD{k}$.
  
 Next we  introduce the formalised operators $ \Rs = \G - \Wa \nWs$ and $\As= K \G + \varrho^2 \Rs$. The operator $\Rs$ is formally equal to $\Pt$ but without boundary restrictions on its domain. The key regularity result we use is the following

\bgP{CompactReg}\hfill

\bgEn{\alph}
\etem If $u \in \WjsD{1} \cap \dom{\Rs}$ then $u \in \R{1}\WlD{1}$.
\etem If $u \in \R{1}\WjsD{0} \cap \dom{\Rs}$ then $u \in \R{2}\WlD{0}$
\enEn
\enP

The operator $\Rs$ is uniformly elliptic on any compact subset of $D$. The proposition follows from careful analysis of the bounding constants of the standard regularity theorems of elliptic operators. See for example, Theorem 8.8 in \cite{Gilbarg:P}. We make the observation that on any sufficiently narrow collar neighbourhood of a compact portion of $\dD$ the function $\varrho$ is equivalent to the Euclidean distance to the boundary. The precise details are tedious but routine, and so will be omitted. 

\bgD{Dirichlet}
For $u,v \in \R{1}\WjsD{0}$ we define the sesquilinear form 
\[ \tQ{v,u}= K \left\{ \G^2  \right\}  \ip{v}{u}{} + \Gm \ip{\varrho v}{\varrho u}{} +  \ip{ \varrho \nWs v}{ \varrho \nWs u}{} + \ip{ (\nW \varrho) v}{ \varrho \nWs u}{} .\]
\enD

\bgL{Basic} 
For sufficiently large $K$ there is a uniform estimate for $u \in \R{1}\WjsD{0}$ of type
\[ \norm[2]{u}{\R{1}\WjsD{0}} \lesssim \left| \tQ{u,u} \right|\]
\enL

For $u \in \CpD$
\begin{align*}
 \ip{\varrho \nWs u}{\varrho \nWs u}{} &= \ip{\nWs (\varrho^2 u) - 2(\nW\varrho)\varrho u}{\nWs u}{}\\
&= \ip{\varrho^2 u}{-\Wa \nWs u}{} - 2\ip{ (\nW\varrho) u}{\varrho \nWs u}{}\\
&= \ip{\varrho^2 u}{(-\nWs\Wa  + [\nWs,W] ) u}{}\\
& \qquad  - 2\ip{ (\nW\varrho) u}{\varrho \nWs u}{}\\
&= \ip{\varrho^2 u}{[(1-\nu)\nWs -\nWs \Ws+ \Ws-\nWs -\uls] u}{} \\
& \qquad - 2\ip{ (\nW\varrho) u}{\varrho \nWs u}{}\\
&=\ip{\varrho^2 u}{(1-\nu-\nWs)\Ws} {} - \ip{\varrho u}{\nu \varrho \nWs u}{}\\
& \qquad  + \ip{\varrho u}{\nu \varrho Wu}{}  - \uls \norm[2]{\varrho u}{}\\
& \qquad - 2\ip{ (\nW\varrho) u}{\varrho \nWs u}{}\\
&= \norm[2]{ \varrho \Ws u}{} + 2\ip{ (W\varrho) u}{\varrho \Ws u}{}  - \nu \ip{ \varrho u}{ \varrho \nWs u}{}\\
& \qquad  + \nu \ip{ \varrho u}{ \varrho \Ws u}{}  - \uls \norm[2]{\varrho u}{}\\
& \qquad - 2\ip{ (\nW\varrho) u}{\varrho \nWs u}{}.\\
\end{align*}
Thus we see that $2\tQ{u,u}$ can be expressed
\begin{align*}
2\tQ{u,u} &=  2K \G^2 \norm[2]{u}{} + (2\Gm  - \uls) \norm[2]{\varrho u}{}  + \norm[2]{\varrho \nWs u}{} + \norm[2]{ \varrho \Ws u}{}\\
& \qquad + 2\ip{ ( W \varrho) u}{\varrho \Ws u}{} - 2\ip{ ( W \varrho) u}{\varrho \nWs u}{} \\
& \qquad  - \nu \ip{ \varrho u}{ \varrho \nWs u}{} + \nu \ip{ \varrho u}{ \varrho \Ws u}{}.  
\end{align*}
We recall from \rfL[TV]{Prelim} that $\uls \lesssim \G^2$ so by a standard small constant, large constant argument we get a bound 
\[ \norm[2]{u}{\R{1}\WjsD{0}} \lesssim \left| \tQ{u,u}\right| \]
for sufficiently large $K$  and $u \in \CpD$. Appealing to a density lemma (\rfT[DY]{Density} ) completes the proof.

\epf

\bgC{WeightedEstimates}
For sufficiently large $K$ if  $u \in \R{1}\WjsD{0} \cap \R{2}\WlD{0} \cap \dom{\As}$ then $u \in \R{2}\WjsD{0}$. Furthermore there exists a constant $C$ independent of $\s \in \Vv$ such that
\[\norm{ u}{\R{2}\WjsD{0}} \lesssim \norm{ \As u}{}\] for all $u \in \R{2}{\WjsD{0}}$.
\enC

\pf
The proof of this result is very similar in nature to that of \rfL[TV]{Estimates} and we shall only sketch out the argument here.

Since $u \in \R{1}\WjsD{0} \cap \dom{\As}$ we see that
\[ \tQ{v,u} = \ip{v}{\As u}{}\] for all $v \in \R{1}\WjsD{0}$. Therefore by \rfL{Basic} we obtain a uniform estimate
\[ \norm{ u}{\R{1}\WjsD{0}} \lesssim \norm{\As u}{}.\]

Let $X$ be any $\s$-uniform operator of order $1$ and $\xi \subset \zeta$ smooth cutoff functions. For any functions $u,v \in \CpD$ we then get
\begin{align*}
 \tQ{v, \varrho X \xi u}&= K \left\{ \G^2  \right\}  \ip{v}{ \varrho X \xi  u}{} + \Gm \ip{\varrho v}{ \varrho^2 X\xi  u}{} +  \ip{ \varrho \nWs v}{ \varrho \nWs \varrho X\xi u}{}\\
 & \qquad  + \ip{ (\nW \varrho) v}{ \varrho \nWs  \varrho X \xi u }{} \\
&= \tQ{ \xi X^* \varrho v,u}  - \Gm \ip{\varrho v}{  \varrho (X\varrho) \xi u}{}+ \ip{\varrho \nWs v + (\nWs \varrho)v}{ [\varrho \nWs, \varrho X \xi ]u}{}\\
& \qquad+ \ip{[\xi X^* \varrho  ,\varrho \nWs] v}{ \varrho \nWs u}{}+ \ip{ \xi (X^*\nWs \varrho) \varrho v}{ \varrho \nWs u}{}\\
& = \ip{X^* \varrho \xi v}{\As \zeta u}{}+ \ip{\varrho L_1 v}{\varrho L_1 \zeta u}{}+ \ip{\varrho L_1 v}{ (\varrho L_1 + \varrho^2 \uls L_0) \zeta u}{}
\end{align*}
Setting $v = \varrho X \xi u$ and applying \rfL{Basic} then implies that
\[ \norm{ \varrho X \xi u}{\R{1}\WjsD{0}} \lesssim \norm{\As \zeta u}{} + \norm{\zeta u}{\R{1}\WjsD{0}} + \norm{\varrho^2 \G^2 \zeta u}{}.\]
The bounding constants can be chosen to depend solely on the functions $\xi,\zeta$ and their derivatives. The final term can be easily controlled in  a manner analogous to \rfL[TV]{Estimates} to yield a uniform estimate
\[ \norm{ \varrho X \xi u}{\R{1}\WjsD{0}} \lesssim \norm{\As \zeta u}{} + \norm{\zeta u}{\R{1}\WjsD{0}} .\]

Next we use the density of \CpD in $\R{2}\WjsD{0}$ (\rfT[DY]{Density}) and the continuity of \As as an operator between $\R{2}\WjsD{0}$ and $L^2(D)$ to see that these are genuine estimates on the whole of $\R{2}\WjsD{0}$. We only a priori know that $u \in \R{2}\WlD{0}$, but since $\xi$ and $\zeta$ have compact support in \hy{2} this is sufficient. The details here mimic \rfT[TV]{Local}.

We can then construct a uniformly locally finite cover of $D$ by hyperbolic balls  as in \rfT[TV]{Transverse} and choose $\xi$ and $\zeta$ to depend on the hyperbolic distance to the centre of each ball. Summing over the cover yields global regularity and the desired estimate.

\epf

\bgT{TangentialW}
If $u \in \dom{\Pt} \cap \{ \Ker{\Pt} \}^\bot$ and $\Pt u \in \WjsD{k}$ then $u \in \R{2}\WjsD{k}$. Furthermore there is a uniform estimate
\[ \norm{u}{\R{2}\WjsD{k}} \lesssim \norm{\Pt u}{\WjsD{k}}.\]
\enT

\pf
Since this result is implied by the stronger estimates of \rfL[TG]{Basic} when $\Gm =0$, we shall suppose throughout that $\Gm \geq \Go >0$.

First suppose $k>0$. Let $X_{k-1}$ be any $\s$-uniform operator of order $k-1$. Then by the work of the previous section we know that $X_{k-1}u \in \WjsD{1}$. Furthermore
\[ \Rs X_{k-1} u  = X_{k-1}\Rs u +(L_k + \ul L_{k-1})u\] we is easily seen to be in $L^2(D)$. Thus by \rfP{CompactReg} we have $X_{k-1} u \in \R{1}\WlD{1}$ and then by \rfC{WeightedEstimates} we see $X_{k-1} u \in \R{1}\WjsD{1}$. This implies $u \in \R{1}\WjsD{k}$ and hence $\varrho u \in \WjsD{k+1}$.

We now recycle this argument letting $X_k$ be any $\s$-uinform operator of order $k$. Then $X_k \varrho u \in \WjsD{1}$ and $\Rs X_k \varrho u$ is easily seen to be in $L^2(D)$. Thus we see $X_k \varrho u \in \R{1}\WjsD{1}$. This is easily sufficient to see that $u \in \R{2}\WjsD{k}$. The uniform estimate is an easy consequence of \rfC{WeightedEstimates} after we make the observation that
\[ \norm{\As u}{\WjsD{k}} \lesssim \norm{\Pt u}{\WjsD{k}}\] by the results of the previous section as $\varrho$ and its derivatives are bounded on $D$.

When $k=0$ we cannot apply the first step of the above argument so must find an alternative way to show that $u \in \R{1}\WjsD{0}$.

Let $\xX$ be the closure of $\CpD$ in $\R{1}\WjsD{0}$ under the norm
\[ \norm[2]{u}{\xX} := \norm[2]{u}{\R{1}\WjsD{0}} + \norm[2]{\nWs u}{}\]
 and consider the sesquilinear form on $\xX$ defined by
\[  \tQ[\xX]{v,u} = \tQ{v,u} + (K + \Gm^{-1})  \ip{\nWs v}{\nWs u}{}\]
By \rfL{Basic} this form is strictly coercive in the sense that there is some $C>0$ such that
\[ \norm{u}{\xX} \leq C \left| \tQ[\xX]{u,u}\right| \]
for all $u \in \xX$. The Lax-Milgram Lemma then implies that if $f \in L^2(D)$ then there exists a unique $u \in \xX$ such that $\tQ[\xX]{v,u} = \ip{v}{f}{}$ for all $v \in \xX$.

For $v \in \xX$ we see that for some constants $C,C\upp$
\begin{align*}
 \left|\ip{\nWs v}{(K\G^2+\varrho^2)\nWs u}{} \right| & \leq C \left| \tQ[\xX]{v,u} \right| + C \left|\tQ{v,u} - \ip{\varrho \nWs v}{\varrho \nWs u}{}\right|\\
 & \leq C\upp \norm{v}{} \left\{ \norm{f}{} + \norm{u}{\xX}\right\}.
 \end{align*}
 We note that $\CpD \subset \xX$ and appeal to \rfL[DY]{Domain} in the appendix to see that $\CpD$ is dense in $\dom{\nWs}$ in the graph norm. Thus the functional $v \mapsto \ip{\nWs v}{(K\G^2 +\varrho^2 )\nWs u}{}$ is bounded on $\dom{\nWs}$ and $(K\G^2+\varrho^2)\nWs u\in \dom{\nWs^*}$. For sufficiently large  $K$  we see  $K\G^2+\varrho^2$ is uniformly bounded below by some positive constant. From the defining properties of $\varrho$ it follows that $(K\G^2+\varrho^2)^{-1}$ is smooth on $\overline{D}$ and has bounded $W$ and $\nW$ derivatives of all orders. It is easy to check that the domain of $\nWs^*$ is closed under multiplication by such functions. Thus we see $\nWs u \in \dom{\nWs^*}$ and $u \in \dom{\Pt}$. The function $u \in \xX \subset \R{1}\WjsD{0}$ is then the unique solution to $\Pt u = \{ K (1+ \Gm^{-1}) +\varrho^2 \}f$. This is sufficient to see that $\dom{\Pt} \subset \R{1}\WjsD{0}$. The rest of the proof now procedes exactly as when $k>0$.

\epf

We can now state a unifying theorem which both yields sharp uniform estimates for  \Pt and addresses the issue of regularity gains in all directions at the boundary and in the interior.

\bgT{Tangential}
Suppose $(D,\nu)$ satisfies either \rfi[TV]{C}{a} or \rfi[TV]{C}{b} and $1 \leq q \leq n-2$. Then for all $\s \in \Vv$ the operators 
\begin{align*} \Pt&: \{ u \in \R{2}\WjsD{k}: \nWs \in \Wo[k+1]\} \cap \{ \Ker{\Pt} \}^\bot \to \WjsD{k} \cap \{ \Ker{\Pt} \}^\bot\\
1&+\Pt: \{ u \in \R{2}\WjsD{k}: \nWs \in \Wo[k+1]\} \to \WjsD{k}
\end{align*}
 are isomorphisms. Furthermore the bounding constants are independent of $\s$.
\enT
The second of these follows from an easy adaption of the techniques we have used to show the first.

\section{Regularity and Estimates for \boxb}\setS{BB}

We shall now restrict our attention to the class of domains described in \rfE[MC]{Model}, i.e. $\Omega = D \times N$ in the strictly pseudoconvex manifold $M=\hy{2} \times N$ . For the  analysis of the preceding sections, this demands that we restrict our attention to precompact sets $D$. We shall permit $\nu$ to be any non-negative constant and it will be suppressed from any notation.

We need weighted versions of the Folland-Stein spaces to deal with  the weighted estimates of the previous section. Accordingly we define
\bgD{Weights}
We define the spaces $\mathscr{V}^{k}\SjO{m}$ inductively by $\mathscr{V}^{0}\SjO{m} = \SjO{m}$ and 
\[ \mathscr{V}^{k+1} \SjO{m} =\left \{ \varphi \in \mathscr{V}^{k}\SjO{m}: \rho \NabH \varphi, \NabH[\top]\varphi \in \mathscr{V}^{k}\SjO{m} \right\}\] with corresponding inductively defined norms.
\enD
Here $\NabH[\top]$ is the tangential component of the connection given explicitly by
\[ \NabH = \NabH[\top] + \nabla_\bt{Y} \otimes \ut^\bt{0} + \nabla_Y \otimes \ut^0.\]
Since we are assuming that $D$ is precompact, we may take $\rho$ to be any smooth defining function for $D$ and identify it naturally with a function on $M$.

There are currently no sharp trace theorems for the Folland-Stein spaces on strictly pseudoconvex CR manifolds. However we can employ the partial Fourier decomposition to develop a useful notion of vanishing at the boundary.
\bgD{Vanishing}
We shall say a $(0,q)$-form   $\varphi \in \SjO{j}$ vanishes on the boundary $\dO$ and write $\varphi \in \So[j]$ if each function $\varphi^\top_\s, \varphi^\bot_\s$ in the partial Fourier decomposition of $\varphi$ is contained in $\Wo[j]$. \enD
It is not difficult to show that $\So[1]$ is equal to the closure of $\Cic{\Omega}$ in $\SjO{1}$. Also $\So[j] =\SjO{j} \cap \So[1]$. We can now describe the domain spaces for our sharp estimate for \boxb.
\bgD{FinalSpace}
\[ \HV{k,j} := \left\{ \varphi  \in \mathscr{V}^{j}\SjO{k}: \varphi^\bot \in \So[k+j], \nabla_{\bt{Y}} \varphi^\top \in \So[k+j-1]  \right\}\]
with norm 
\[ \norm[2]{\varphi}{\HV{k,j}} := \norm[2]{\varphi}{\mathscr{V}^{j}\SjO{k}} + \norm[2]{\varphi^\bot}{\SjO{k+j}} + \norm[2]{\nabla_{\bt{Y}} \varphi^\top}{\SjO{k+j-1}}.\]
\enD
In particular, we note that $\SjO{k+j} \subset \HV{k,j} \subset \SjO{k}$ with the inclusion maps continuous.

Now set  $\Kr:= \Ker {\boxb}$ with $\boxb$ regarded as an unbounded operator on $L^2(\Lambda_\ut^{0,q}\Omega)$. Then since \[\ip{\boxb \varphi}{\varphi}{L^2(\Omega)} = \ \norm[2]{\db \varphi}{L^2(\Omega)} + \norm[2]{\dbs \varphi}{L^2(\Omega)}\] it follows immediately that $\Kr = \Ker{\db} \cap \Ker{\dbs}$.

\bgT{Main}
Suppose $\Omega$ is as in \rfE[MC]{Model}. Fix  $1 \leq q \leq n-2$.  Then for each $k \geq0$, the operators
\begin{align*}
 \boxb&: \HV{k,2} \cap (\Kr)^\bot \to \SjO{k} \cap (\Kr)^\bot\\
 1&+ \boxb: \HV{k,2} \to \SjO{k}
\end{align*}
 are both isomorphisms on $(0,q)$-forms.  
\enT

\pf Expand out the $(0,q)$-form  $\varphi$ in terms of the partial Fourier decomposition and decompose the relevant norms using \rfL[MC]{FourierFS} . Apply \rfT[TV]{Transverse} to each of the tranverse components and \rfT[WR]{Tangential} to each of the tangential components of $\varphi$. The constants are uniform across all choices of $\s$ so a simple comparison theorem argument completes the proof.
\epf

\bgC{NonHomogeneous}
Suppose the $(0,q)$-form $\vs \in \SjO{k}$  satisfies $\db \vs =0$ and $\vs \bot \Kr$. Then there exists a unique $\varphi \in \HV{k,1}$ such that $\varphi \bot \Kr$ and $\db \varphi = \vs$. Furthermore there is a constant $C$ independent of $\vs$ such that \[\norm{\varphi}{\HV{k,1}} \leq C\norm{\vs}{\SjO{k}}.\]
\enC

\pf
By \rfT{Main} there exists a unique $\ua \in \HV{k,2} \bot \Kr$ such that $\boxb \ua = (\db \dbs + \dbs \db) \ua =\vs$ and $\norm{\ua}{\HV{k,2}} \leq C\upp \norm{\vs}{\SjO{k}}$ for some $C\upp$ independent of $\vs$. From $\db \vs =0$ it follows that $\vs \bot \rng{\dbs}$ and so $\vs = \db \dbs \ua$. Take $\varphi = \dbs \ua$. From the proof of \rfL[MC]{KohnLaplacian} we know that $\dbs[\top]$ and $\nabla_\bt{Y}$ commute. Tangential derivatives map $\SjO{k+1}$ to $\SjO{k}$. Thus from \rfL[MC]{Facts} it follows that $\dbs$ is a bounded linear operator from $\HV{k,2}$ to $\HV{k,1}$.

\epf

\bgL{Cohomology}
Suppose $1\leq q \leq n-2$. If  $\hKR{0,q}{N}=0$ then $\Kr = 0$ otherwise $\dim \Kr = \infty$.
\enL

\pf
A $(0,q)$-form $\varphi$ is in $\Kr$ if and only if $\varphi^\bot=0$ and $\varphi^\top$ expands as a sum of forms $\varphi^\top_\s \s$ such that $\Gm=0$ and $\nWs \varphi^\top_\s =0$. Now $\{\s \in \Vv: \Gm =0\} \cong \hKR{0,q}{N}$. If this set is non-empty  then we can construct elements of $\Kr$ from any holomorphic function multiplied by an appropriate power of $s$.  
\epf

\bgC{Cohomology}
If $\hKR{0,q}{N}=0$ then $\hKR{0,q}{\Omega} = 0$ otherwise $\dim \hKR{0,q}{\Omega} = \infty$.

\enC

\pf
This now follows immediately from the observation that $\Kr \bot\rng{\db}$.
\epf

\bgL{Hypoelliptic}
The  operator $\boxb$ is hypoelliptic up to the boundary  if and only if $\hKR{0,q}{N}=0$.
\enL

\pf
If $\hKR{0,q}{N}=0$ then $\boxb$ is injective by \rfL{Cohomology}. Thus by \rfT{Main} if $\boxb \varphi \in \Ci{\Omega}$ then $\varphi = \bigcup\limits_{k \geq 0} \SjO{k}$. Since there is a continuous inclusion $\SjO{2k} \subset \HjO{k}$, hypoellipticity follows from the Sobolev Lemma.

If $\hKR{0,q}{N} \ne 0$ then we can choose some $\s \in \Vv$ such that $\Gm =0$, then following \rfX[TG]{NegRegularity} we can easily construct $\varphi \in \Kr$ such that $\varphi \notin \SjO{1}$.

\epf

It is worth pointing out here that regardless of the cohomology of the foliating manifold $N$, the operator $1+\boxb$ is always hypoelliptic up to the boundary of $\Omega$. Likewise $\boxb$ is always hypoelliptic in the interior.

\bgL{NonCompact}
For any $1\leq q \leq n-2$ The bounded operator $(1+\boxb)^{-1}: L^2(\Omega) \to L^2(\Omega)$ is not compact on $(0,q)$-forms.
\enL

\pf If $\hKR{0,q}{N}\ne 0$ then this can be seen from \rfL{Cohomology}. However in general  
 this is really a statement about the operators $\Pt$. Fix $w_0 \in D$ and $\e>0$ such that the disc $|w-w_0|<\e$ is contained in $D$. Choose any $\s \in \Vv$  and let $f$ be a function on $D$ that can be expressed as the product of $s^{(\uls+\nu)/2}$ and a holomorphic function. Then $(1+\Pt) f = (1+\Gm) f$. Now consider the sequence $f_k$ where $f_k= s^{(\uls+\nu)/2} c_k w^k$ where $c_k$ is chosen to make $\norm{f_k}{L^2(D)}=1$. If some subsequence $f_j$ converges in $L^2(D)$ then the same subsequence must converge in $L^2(|w-w_0|<\e)$. However the functions $(w-w_0)^j$ are mutually orthogonal on $|w-w_0|<\e$ in the Euclidean $L^2$ metric. Since $s$ is bounded above and below on this disc,  it follows that the subsequence $f_j$ cannot converge on $|w-w_0|$ in the hyperbolic $L^2$ metric. Thus the bounded sequence $f_k$ has no convergent subsequences and so  $(1+\Pt)^{-1}$ is not compact as an operator on $L^2(D)$. 
 
\epf

As an immediate consequence of this, we see
\bgC{NonSub}
The operator $\boxb$ is not globally subelliptic on $\Omega$, i.e. there is no $\e >0$ such that
\[ \norm{\varphi}{\HjO{\e}} \leq C \norm{\boxb \varphi}{L^2(\Omega)} + \norm{\varphi}{L^2(\Omega)}.\]
\enC

For completeness, in addition to the global theory, we can combine the results of \rfT[TV]{Local} and \rfC[TG]{Local} to establish the following local estimates for $\boxb$.

\bgT{Local}
Under the same conditions as \rfT{Main}, then for all $k \geq 0$ and any pair of smooth cut-off  functions  $\xi \subset \zeta$ on $M$ that depend solely on $w$ there is a constant $C>0$ such that
\[ \norm{ \xi \varphi}{\HV{k,2}} \leq C \norm{\zeta}{\SjO{k}} + \norm{\zeta \varphi}{\lO} + \norm{ \nabla_\bt{Y} \varphi^\top}{\lO}\] whenever $\varphi \in \HV{k,2}$. 
\enT

\appendix
 \section{Density Results}\setS{DY}
 
This section is dedicated to proving the various density results used at various stages throughout this paper. We shall often use $\Wj{k}$ when a statement holds for all $\WjsD{k}$. We shall use the notation $\R{k}\Wjs[]{j}$ to describe weighted spaces exactly as in \rfD[WR]{Weighted}, but shall only insist that the weighting function $\rho$ be a  smooth defining function for $D$.
 
 \bgD{Snorm}
 The pointwise norm associated to $\Wj[]{m}$ of a sufficiently differentiable function is given by 
\[ \snorm[2]{f}{\Wj[]{m}} = \sum\limits_{j+k \leq m} \snorm[2]{ W^j \nW^k f}{}.\]
\enD

Throughout this section let $D$ be a smooth open domain in \hy{2} such that the functions $s$ and $t$ are bounded on $D$. We define $ D_{-}^\e = \{w \in D:s < \e \}$ and $D_+^\e  = \{w \in D:s > \e \}$. 

\bgD{Plus}
If $\xX$ is a space of functions defined on $D$ then $\xX_+$ is defined as
\[ \xX_+ = \{f \in \xX: \text{there is some $\e>0$ such that $f=0$ on $D_{-}^\e$}\}.\]
\enD

\bgL{Plus}
There exists a family of smooth functions $\{\xi_n\}$ in \CpD  and constants  $C_1,C_2, \dots$ such that 
\bgEn[L]{\arabic}
\etem Each $\xi_n$ is supported in $D^{1/n}_+$ and is identically one on $D^{2/n}_+$
\etem $\snorm{\xi_n}{\Wj[]{k}} \leq C_k$  for all $n$.
\enEn
\end{lemma}

\pf
 Fix a smooth function $\xi$ on the real line such that $|\xi(x)| \leq 1$ everywhere, $\xi(x)=0$ on $|x|\leq 1$ and $\xi(x)=1$ on $|x|\geq 2$. Set $c_k$ to be an upper bound for $\xi^{(k)}$ on \rn{}. Now define $\xi_n(w) = \xi(ns)$. Then 
\[
\left| (\nW^k \xi_n )(w)\right| \leq 2^k \sum_{j \leq k} \left| (sn)^j \xi^{(j)}(ns) \right|.\]
But $\xi^{(j)}(ns)=0$ unless $1 < ns < 2$. So we see
\[
\left| (\nW^k \xi_n )(w)\right| \leq 2^k \sum_{j \leq k} 2^j c_j.\]
The derivatives of other types can be controlled in a similar fashion.

\epf

\bgC{Plus}\hfill

\bgEn{\alph}
\etem $\Wpj[D]{k}$ is dense in $\Wj[D]{k}$. 
\etem If $\sup\limits_D \snorm{\rho}{\Wj{k}} < \infty$ then $\Rp{k}L^2(D)$ is dense in $\R[]{k}L^2(D)$.
\enEn
\end{cor}

\pf
In either case it is easy to see that $\xi_n f \to f$ in the relevant norm.

\epf

\bgT{Density}
If $\sup\limits_D \snorm{\rho}{\Wj{k}} < \infty$ then \CpD is dense in $\R[]{k}L^2(D)$.
\end{thm}

\pf The proof runs by induction. The condition that $\sup \snorm{\rho}{ \Wj{k}} < \infty$ is necessary only to guarantee that $\Wj{k} \subset \R[]{k}L^2(D)$. By \rfC{Plus} it is sufficient to show $\CpD$ is dense in $\Rp{k}L^2(D)$. 

Let $N$ be the gradient of $\rho$ viewed purely as a function of $w$. Then $N$ is a smooth vector field that is non-vanishing in a region $ \{-\e<\rho < \e\}$ for some $\e >0$. Let
$\Theta$ be the maximal flow for $\xi N$ where $\xi$ is a smooth function with $0 \leq
\xi \leq 1$, $\xi =1$ on $|\rho| < \frac{\e}{2}$ and
$\xi=0$ on $|\rho| \geq \e$. For a function $f$ on $D$ and small
$\ud>0$,  define a new
function $f^\ud$ by $f^\ud(z) = f \circ \Theta(z, \ud)=
\Theta_{\ud }^* f$. 

\hfill

\noindent \textbf{Claim:} If $u \in \Rp{k}L^2(D)$ then $u^\ud \in \Wj{k}$ and $u^\ud \to u$ in $\R[]{k}L^2(D)$.

\hfill

\noindent Sub-claim 1: For sufficiently small $\ud_0$ and $\e>0$, the operators $\Theta_{\ud}^*$ for $0<\ud<\ud_0$ are uniformly bounded on $L^2(D^\e_+)$.

\hfill

Recall that $dV = s^{-2} ds \wedge dt$. Since $\Theta_{- \ud}$ is
  orientation-preserving, there is a smooth function $B_\ud$ such that
  $B_\ud \geq 0$ and $\Theta_{-  \ud}^* dV = B_\ud dV$. The functions $B_\ud$ are smooth everywhere and converge uniformly to 1 on the compact set $D_+^\e$. Thus for some $\ud_0$, the functions $B_d$ are uniformly bounded above on $D_+^\e$ by some constant $C_0$ for all $0<\ud<\ud_0$.

\begin{align*}
 \norm[2]{f^\ud}{L^2(D^\e_+)} &= \int\limits_{D^\e_+} \snorm[2]{\Theta_{ \ud}^* f}{} \, dV = 
\int\limits_{\Theta_\ud(D^\e_+)} \Theta_{- \ud}^* \big( \Theta_{ \ud} \snormV[2]{f}{} \, dV \big)\\
& =\int_{\Theta_\ud(D^\e_+)} \snormV[2]{f}{} \Theta_{- \ud}^* dV = \int_{\Theta_\ud(D^\e_+)} B_\ud \snormV[2]{f}{} \, dV\\
& \leq C_0 \normV[2]{f}{L^2(\Theta_\ud(D^\e_+))} \leq C_0 \normV[2]{f}{L^2(D^\e_+)}.
\end{align*}
\espf

This allows us to prove that the claim is true for $k=0$. Choose $u$ in $\Rp{k}L^2$ and fix $\e$ such that $u$ is supported in $D^{2\e}_+$. Thus for sufficiently small $\ud$, $u^\ud$ is supported in $\{s>\e\}$. If $h \in \Cic{D}$ then $h$ is uniformly continuous and $h^\ud$ converges uniformly to $h$. We then recall that any $L^2$ function on a compact domain can be approximated by compactly supported smooth functions. Finally apply sub-claim (1) to observe that $\|{h^\ud -g^\ud}\|_{L^2(D)} \leq C_0 \norm{h-g}{L^2(D)}$. A simple triangle inequality argument completes this step.

Now suppose the main claim is true for all $j<k$ and fix $u \in \Rp{1}L^2(D)$. Choose an $\e$ such that $u$ vanishes on $D^{2\e}_-$.

\hfill

\noindent Sub-claim 2: For small $\ud$, $u^\ud \in \Wj{k}$ and $(Wu)^\ud$, $(\nW u)^\ud \in \Wj{k-1}$.

\hfill
We argue by direct computation. The vector field $\Theta_{\ud*}W$ is a smooth combination of $W$ and $\nW$. Thus
\begin{align*}
\int\limits_D \left| W^k u^\ud \right|^2 dV &= \int\limits_{\Theta_\ud(D)} \Theta_{-\ud}^* \left| W^k u^\ud \right|^2 B_\ud dV\\
&= \int\limits_{\Theta_\ud(D)}  \left| (\Theta_{\ud *} W)^k u \right|^2 B_\ud dV\\
&\leq C_0 \int\limits_{\Theta_\ud(D^\e_+)}  \left| (\Theta_{\ud *} W)^k u \right|^2  dV < \infty.
\end{align*} 
where the final inequality follows from the observation that $u \in \Wj[\Theta_\ud(D^\e_+)]{k}$. The other derivatives and cases are proved similarly.

\espf

\noindent Sub-claim 3: For small $\ud$ there exists a constant $C_\rho$ such
that $\sup\limits_{D^\e_+} |\rho^\ud -\rho|<C_\rho \ud$.

\hfill

Since $\Theta$ is the flow of the gradient of $\rho$ near $\partial D$, we have 
\begin{align*}
0\leq \rho^\ud(z) - \rho(z) &= \int_0^{ \ud} \aip{\xi N (\Theta(z,a))}{
\frac{d}{da} \Theta(z,a)}{} \, da\\
& = \int_0^{ \ud}   |\xi N |^2 da\\
& \leq \left(\sup_{D^\e_+} |\xi N|^2\right)  \ud 
\end{align*} 
\espf

\noindent Sub-claim 4: $\rho (Wu)^\ud - \rho W u^\ud \to 0$ in $\R[]{k-1}L^2(D)$. A similar result holds for $\nW$.

\hfill

Again we compute directly that
\begin{align*}
\frac{1}{\ud^2} \int\limits_D \left| (\rho W)^{k-1} \rho ((Wu)^\ud - Wu^\ud)  \right|^2 dV &= \int\limits_{\Theta_\ud(D)} \Theta_{-\ud}^* \left| (\rho W)^{k-1} \rho \frac{1}{\ud} \Theta_\ud^* \left( (W-\Theta_{\ud *}W)u \right) \right|^2 B_\ud dV\\
&= \int\limits_{\Theta_\ud(D)}  \left| (\Theta_{\ud *} \rho W)^k \left(\rho^{-\ud} \frac{1}{\ud} (W-\Theta_{\ud *}W) u\right) \right|^2 B_\ud dV.
\end{align*} 
Now as $W$ is a vector field with smooth coefficients we see that $\Theta_{\ud *} \rho W \to \rho W$ and $ \frac{1}{\ud} (W-\Theta_{\ud *}W)  \to [N,W]$, both locally uniformly. In addition $u$ is supported in $D^\e_+$. Therefore this integral can be uniformly bounded as $\ud \to 0$. Thus \[\norm{\rho ((Wu)^\ud - Wu^\ud)}{\R[]{k-1}L^2(D))} \longrightarrow 0.\]
The other derivatives are then handled in a similar fashion.

\espf

We are now ready to prove the main claim. The first part of the statement was proved in sub-claim (2) so it only remains to show the convergence property. But
\begin{align*}
\norm{\rho W (u^\ud -u)}{\R[]{k-1}L^2} &\leq \norm{ (\rho W u)^\ud - (\rho W u)}{\R[]{k-1}L^2} + \norm{ \rho (Wu)^\ud - \rho Wu^\ud}{\R[]{k-1}L^2}\\
& \qquad + \norm{ (\rho^\ud -\rho)(Wu)^\ud}{\R[]{k-1}L^2}.
\end{align*}
The first term tends to zero by induction, the second by sub-claim (4) and the third by sub-claims (2) and (3).

\espf

To complete the proof of the theorem itself we note that  each  $\Wj{k}$ is equivalent to the standard hyperbolic Sobolev space of order $k$ on $D$. Thus for $\ud$ sufficiently small, $u^\ud$ can be approximated in the $\Wj{k}$ norm by functions in $\CiD$. Since $u^\ud$ is easily seen to be  in $\Rp{k}L^2(D)$ we can replace \CiD with \CpD.  The proof is then completed by a straightforward diagonalisation argument.

\epf

\bgL{Domain}
With $\nW$ considered as an unbounded operator $L^2(D)$ to $L^2(D)$, $\CpD$ is dense in $\dom{\nW}$ in the graph norm $\| \cdot\| + \|\nW \cdot\|$.
\end{lemma}

\pf By \rfC{Plus} it is sufficient to prove the result for $\domp{\nW}$. The operator $\nW$ is elliptic so $\dom{\nW} \subset \WjMain{1}{loc}{(D)}$. Therefore if $U$ is any precompact open set with $\overline{U} \subset D$ then $u_{|U} \in \Wj[U]{1}$ whenever $u \in \dom{\nW}.$

Thus following the proof of \rfT{Density} we see that if $u \in \domp{\nW}$ then $u^\ud \in \Wj{1}$. The argument for sub-claim 4 in the previous lemma can easily be adapted to complete the proof.

\epf

\bgD{CpO}For fixed $0 \leq q \leq n-1$ we define the space $\CpO$ to be the collection of $(0,q)$-forms in $\lO$ such that each function of the partial Fourier decomposition is in $\CpD$.
\end{defn}

\bgC{CpOdensity}
The space $\CpO$ is dense in $\dom{\db}$ in the graph norm.
\end{cor}

\pf
If a $(0,q)$-form $\varphi$ is in $\dom{\db}$ then each element of its partial Fourier decomposition is in $\dom{\nW}$. Thus it can be approximated arbitrarily closely by a smooth function on $\overline{D}$. A simple argument then allows us to construct the partial Fourier decomposition of a form in \CpO which approximates $\varphi$ arbitrarily closely in the graph norm of $\db$.
 
\epf

\bgL{Vanish}
The set $\Cv{D}=\{u \in \CpD: u=0 \text{ on $\partial D$}\}$ is dense in $\Won{k}{D}$.\end{lemma}

\pf Again it is sufficient to show that any  $\WoX{k}{(D)}{+}$ can be approximated by compactly supported smooth functions. Any such $u$ vanishes on some $D^{2\e}_-$ it can be viewed as an element of $\Wj[D\upp]{k}$ where $D\upp$ is a smoothly bounded precompact domain such that the part of the boundary above $s=\e$ coincides with the boundary of of $D$ there. This new domain $D\upp$ can easily be chosen so that $u \in \Won{1}{D\upp}$.

As $D\upp$ is compact it admits a well-defined continuous trace function from $\Wj[D\upp]{k}$ to $\Wj[\partial D\upp]{k-1/2}$ with continuous right inverse. Let $v_n \in \Ci{D\upp}$ converge to $u \in \Wj[D\upp]{k}$. As $u \in \Won{1}{D\upp}$ the trace of $u$ is identically zero. As the trace is continuous this implies that the traces of the functions $v_n$ converge to zero in $\Wj[\partial D\upp]{k-1/2}$. The right-inverse of the trace maps smooth functions to smooth functions and is continuous into the $\Wj[D\upp]{k}$ norm. Thus when we extend the traces of the sequence $v_n$ back onto $D\upp$ the resulting functions converge to zero in the $\Wj[D\upp]{k}$ norm. Subtracting these functions from $v_n$ then produces a sequence of functions in \Cv{D\upp} that converge to $u$ in $\Wj[D\upp]{k}$. But since $u$ vanishes on $\{s<2\e\}$ these smooth functions can easily be extended to $\Cv{D}$ in such a way to obtain convergence in $\Wj{k}$.

\epf

\bibliographystyle{plain}
\bibliography{References}

\end{document}